\documentclass[leqno,11pt]{amsart}

\usepackage[top=3.75cm, bottom=3cm, left=3cm, right=3cm]{geometry}
\usepackage[dvipsnames, table]{xcolor}
\usepackage{amsthm}
\theoremstyle{plain}
\usepackage{amssymb}
\usepackage{marvosym}
\usepackage{bm}
\usepackage{mathrsfs}
\usepackage{enumerate}  
\usepackage{mathtools}
\usepackage{tikz-cd}
\usepackage{graphicx}
\usepackage{float}
\usepackage[titletoc,title]{appendix}
\usepackage{marginnote}
\usepackage[disable]{todonotes}
\usepackage{etoolbox}
\usepackage[all]{xy}
\makeatletter
\patchcmd{\Ginclude@eps}{"#1"}{#1}{}{}
\makeatother
\usepackage[outdir=./]{epstopdf}
\usepackage[numbers]{natbib}
\usepackage[utf8]{inputenc}
\usepackage[english]{babel}
\usepackage{changepage}
\usepackage{makecell}
\usepackage{enumitem}
\usetikzlibrary{decorations.pathmorphing}
\usepackage{cancel}

\definecolor{lightblue}{HTML}{1F88CD}
\definecolor{lightgrey}{HTML}{727272}
\definecolor{lightblue2}{HTML}{009EC1}
\definecolor{mypink}{HTML}{FD00B0}
\definecolor{lightred}{HTML}{ff4d4d}
\def\A{\mathcal{A}}

\usepackage{hyperref}
\hypersetup{
	colorlinks=true,
    linkcolor={OliveGreen},
    citecolor={lightred},
	urlcolor={black}
}

\newtheorem*{theorem*}{Theorem}
\newtheorem{theorem}{Theorem}[section]
\newtheorem{corollary}[theorem]{Corollary}
\newtheorem{lemma}[theorem]{Lemma}

\newtheorem{proposition}[theorem]{Proposition}
\theoremstyle{definition}

\theoremstyle{definition}
\newtheorem{definition}[theorem]{Definition}
\theoremstyle{definition}
\newtheorem{remark}[theorem]{Remark}
\theoremstyle{definition}

\theoremstyle{definition}

\theoremstyle{definition}

\theoremstyle{definition}

\theoremstyle{definition}

\theoremstyle{definition}
\newtheorem{question!}[theorem]{Question!}
\theoremstyle{definition}

\makeatletter
\newcommand*\sbt{\mathpalette\sbt@{.75}}
\newcommand*\sbt@[2]{\mathbin{\vcenter{\hbox{\scalebox{#2}{$\m@th#1\bullet$}}}}}
\makeatother

\newcommand{\D}{\mathrm{D}}

\newcommand{\ZZ}{\mathbb{Z}}

\newcommand{\PP}{\mathbb{P}}

\newcommand{\ch}{\mathrm{ch}}

\newcommand{\HH}{\mathrm{HH}}

\DeclareMathOperator{\im}{im}

\DeclareMathOperator{\Hom}{Hom}

\DeclareMathOperator{\Pic}{Pic}

\DeclareMathOperator{\HP}{\mathrm{HP}}
\newcommand{\bv}{\mathbf{v}}
\newcommand{\cA}{\mathcal{A}}

\newcommand{\Ku}{\mathcal{K}u}

\DeclareMathOperator{\oh}{\mathcal{O}}

\title[Categorical Torelli theorem via equivariant Kuznetsov component]{Three approaches to a categorical Torelli theorem for cubic threefolds of non-Eckardt type via the equivariant Kuznetsov components}  

\author{Sebastian Casalaina-Martin, Xianyu Hu, Xun Lin, Shizhuo Zhang and Zheng Zhang}

\begin{document}

\begin{abstract}
	Let $Y$ be a cubic threefold with a non-Eckardt type involution $\tau$. Our first main result is that the $\tau$-equivariant category of the Kuznetsov component $\Ku_{\mathbb{Z}_2}(Y)$ determines the isomorphism class of $Y$ for general $(Y,\tau)$. We shall prove this categorical Torelli theorem via three approaches: a noncommutative Hodge theoretical one (using a generalization of the intermediate Jacobian construction in \cite{perry2020integral}), a Bridgeland moduli theoretical one (using  equivariant stability conditions), and a Chow theoretical one (using some techniques in \cite{kuznetsovnonclodedfield2021}).The remaining part of the paper is devoted to proving an equivariant infinitesimal categorical Torelli for non-Eckardt cubic threefolds $(Y,\tau)$. To accomplish it, we prove a compatibility theorem on the algebra structures of the Hochschild cohomology of the bounded derived category $D^b(X)$ of a smooth projective variety $X$ and on the Hochschild cohomology of a semi-orthogonal component of $D^b(X)$. Another key ingredient is a generalization of a result in \cite{macri2009infinitesimal} which shows that the twisted Hochschild-Kostant-Rosenberg isomorphism is compatible with the actions on the Hochschild cohomology and on the singular cohomology induced by an automorphism of $X$. In appendix, we prove an equivariant categorical Torelli theorem for arbitrary cubic threefold with a geometric involution under a natural assumption.   
\end{abstract}

\subjclass[2010]{Primary 14F05; secondary 14J45, 14D20, 14D23}
\keywords{Derived categories, Kuznetsov components, categorical Torelli theorems, cubic threefolds, equivariant categories, equivariant stability conditions.}

\address{Department of Mathematics, University of Colorado, Boulder, CO 80309, USA}
\email{casa@math.colorado.edu}

\address{Fakult\"at f\"ur Mathematik,
Technische Universit\"at M\"unchen, D-85747 Garching bei M\"unchen, Germany}
\email{xianyu.hu@tum.de}

\address{Max Planck Institute for Mathematics, Vivatsgasse 7, 53111 Bonn, Germany}
\email{xlin@mpim-bonn.mpg.de}

\address{Simons Laufer Mathematical Sciences Institute,17 Gauss Way, Berkeley, CA 94720}
\email{shizhuozhang@msri.org, shizhuozhang@mpim-bonn.mpg.de}
\address{Institut de Mathématiqes de Toulouse, UMR 5219, Université de Toulouse, Université Paul Sabatier, 118 route de
Narbonne, 31062 Toulouse Cedex 9, France}
\email{shizhuo.zhang@math.univ-toulouse.fr}

\address{Institute of Mathematical Sciences, ShanghaiTech University, Shanghai 201210, China}
\email{zhangzheng@shanghaitech.edu.cn}
\date{\today}

\maketitle

{
\hypersetup{linkcolor=blue}
\setcounter{tocdepth}{1}
\tableofcontents
}

\section*{Introduction}
	Torelli problems are some of the oldest and the most classical problems in algebraic geometry. The classical Torelli problem asks whether an algebraic variety $X$ is uniquely determined by the Hodge structure on its cohomology. Specifically, denote by $\mathcal{P}$ a period map $\mathcal{P}:\mathcal{M}\rightarrow\mathcal{D}/\Gamma$, where $\mathcal{M}$ is a moduli space parameterizing isomorphism classes of certain algebraic varieties $X$, and $\mathcal{D}/\Gamma$ is a period domain (that is, a classifying space of some polarized Hodge structures). One says that a \emph{global Torelli theorem} holds if $\mathcal{P}$ is injective, and that an \emph{infinitesimal Torelli theorem} holds if the differential $d\mathcal{P}$ of $\mathcal{P}$ is injective. Torelli theorems hold for wide class of varieties, see \cite{catanese1984torelli}. For example, the period map for smooth cubic threefolds $Y\subset \PP^4$ is defined by sending $Y$ to its intermediate Jacobian $J(Y)$ (or equivalently, the weight $1$ polarized Hodge structure on $H^3(Y,\ZZ)(1)$); both global and infinitesimal Torelli theorems hold for smooth cubic threefolds.

	A categorical variant of the classical Torelli problem, also referred to as a \emph{categorical Torelli problem}, asks whether the non-trivial semi-orthogonal component $\Ku(X)$ (called the \emph{Kuznetsov component}) of the bounded derived category $D^b(X)$ of an algebraic variety $X$ determines the isomorphism class of $X$. Similar to classical Torelli theorems, it has been shown that categorical Torelli theorems hold for many interesting varieties including Enrqiues surfaces (\cite{li2021refined}, \cite{li2022refined}), cubic threefolds and cubic fourfolds (\cite{bernardara2012categorical}, \cite{huybrechts2016hochschild}, \cite{feyzbakhsh2023new}), (weighted) Fano hypersurfaces (\cite{pirozhkov2022categorical}, \cite{Lahoz2023categorical}, \cite{Rennemo2023hochschild}, \cite{lin2023serre}) and several prime Fano threefolds (\cite{jacovskis2021categorical}, \cite{dell2023cyclic}). We refer the reader to \cite{PS2023cattorelli} for more details. 

	In the present article, we consider categorical Torelli problems for algebraic varieties with additional automorphisms. More specifically, we focus on cubic threefolds $Y$ admitting a \emph{non-Eckardt type} involution $\tau$ which is a biregular involution with fixed locus the disjoint union of a cubic curve $C$ and a line $L$. Cubic threefolds $(Y,\tau)$ with a non-Eckardt type involution have been studied in \cite{casalaina2022moduli}; among their results, the authors show that $(Y, \tau)$ is determined up to isomorphism by the invariant part $J(Y)^{\tau}$ of the intermediate Jacobian. From the categorical perspective, we consider the semi-orthogonal decomposition of the bounded derived category $D^b(Y)$
	$$D^b(Y)=\langle\Ku(Y),\oh_Y,\oh_Y(1)\rangle$$
where $\Ku(Y)$ denotes the Kuznetsov component, the orthogonal complement of the line bundles $\oh_Y$ and $\oh_Y(1)$. The geometric involution $\tau$ induces an action of the group $\mathbb{Z}_2=\langle 1,\tau\rangle$ on $D^b(Y)$ preserving $\Ku(Y)$, and we denote the corresponding equivariant Kuznestov category by $\Ku_{\mathbb{Z}_2}(Y)$. As $\Ku(Y)$ determines the isomorphism class of $Y$, it is a natural question whether the equivariant Kuznetsov component $\Ku_{\mathbb{Z}_2}(Y)$ determines the isomorphism class of the non-Eckardt type cubic threefold $(Y,\tau)$. The main result of our paper in this direction is to give an affirmative answer to this question for general cubic threefolds with a non-Eckardt type involution. 
 
\begin{theorem}[=Theorems \ref{theorem_equivariant_categorical_Torelli}, \ref{alterntative_proof}]\label{main_theorem_first_equivariant_categorical_Torelli}
	Let $(Y,\tau)$ and $(Y',\tau')$ be general cubic threefolds with a non-Eckardt type involution. Assume that there is a Fourier-Mukai type equivalence $\Phi:\Ku_{\mathbb{Z}_2}(Y)\simeq\Ku_{\mathbb{Z}_2}(Y')$ between the equivariant Kuznetsov components. Then $(Y,\tau)\cong (Y',\tau')$.
\end{theorem}

	We shall study the equivariant Kuznetsov components for non-Eckardt cubic threefolds from three different perspectives: a noncommutative Hodge theoretical one (using an equivariant version of the intermediate Jacobian construction in \cite{perry2020integral}), a moduli space theoretical one (via Bridgeland moduli spaces of stable objects in the equivariant Kuznetsov components), and a Chow theoretical one (using Abel--Jacobi maps and some techniques in \cite{kuznetsovnonclodedfield2021}). Consequently, we prove Theorem \ref{main_theorem_first_equivariant_categorical_Torelli} (however, we are only able to prove a slightly weaker version of Theorem \ref{main_theorem_first_equivariant_categorical_Torelli} in the Chow approach, see Theorem \ref{alterntative_proof_chow}). More details are given in Subsections \ref{section_K_theory}, \ref{section_Bridgeland} and \ref{section_Chow}. 

	The upshot is that it looks to us that using the equivariant Kuznetsov component $\Ku_{\mathbb{Z}_2}(Y)$, from the first two perspectives one can only determine the invariant part $J(Y)^{\tau}$ of the intermediate Jacobian for a general $Y$, while from the third perspective, one can only determine $J(Y)^{\tau}$ up to isogeny (which necessitates the assumption, in that approach, that $(Y,\tau)$ be \emph{very general}). In a little more detail, which we expand on in the subsections below, the situation is as follows. As mentioned above, from the results of \cite{casalaina2022moduli}, one wants to show that the equivariant Kuznetsov component $\Ku_{\mathbb{Z}_2}(Y)$ naturally determines the invariant part of the intermediate Jacobian, $J(Y)^\tau$. From the noncommutative Hodge theoretic perspective, and the moduli space theoretic perspective, one finds that while the Kuznetsov component $\Ku(Y)$ naturally determines the Fano variety of lines on $Y$, and consequently the intermediate Jacobian $J(Y)$, the equivariant Kuznetsov component $\Ku_{\mathbb{Z}_2}(Y)$ naturally determines the invariant lines on $Y$, which are parameterized by a genus $4$ bielliptic curve $\widetilde C$, together with an isolated point. There is a surjection $J(\widetilde C)\twoheadrightarrow J(Y)^\tau$, with kernel an elliptic curve (this can be seen from the Prym construction for the cubic threefold), but unfortunately, we were not able to recover this surjection given only the data of $\Ku_{\mathbb{Z}_2}(Y)$. For $Y$ general, we can recover $J(Y)^\tau$ up to isomorphism, however, from $J(\widetilde C)$. In a different direction, there exists a genus $4$ bielliptic curve carrying more than one bielliptic structure (cf.~\cite{CDC05}), and hence by \cite{casalaina2022moduli} there exist non isomorphic non-Eckardt cubic threefolds $(Y,\tau)$ and $(Y',\tau')$ with $J(\Ku_{\mathbb{Z}_2}(Y))\cong J(\Ku_{\mathbb{Z}_2}(Y'))$. However, it is not clear to us at this moment if $\Ku_{\mathbb{Z}_2}(Y)\simeq \Ku_{\mathbb{Z}_2}(Y')$ and whether one can remove the assumption that non-Eckardt cubic threefolds are general in Theorem \ref{main_theorem_first_equivariant_categorical_Torelli}. From the Chow theoretical perspective, while $\Ku(Y)$ naturally determines the integral Hodge structure $J(Y)$, we are only able to show that the equivariant Kuznetsov component $\Ku_{\mathbb{Z}_2}(Y)$ naturally determines the $\mathbb Z[1/2]$-Hodge structure $J(Y)^\tau_{\mathbb Z[1/2]}$; this allows us to prove Theorem \ref{main_theorem_first_equivariant_categorical_Torelli} for very general cubic threefolds carrying a non-Ekcardt type involution. 

\subsection{Topological K-theory and Hodge theory for the equivariant Kuznetsov components}\label{section_K_theory} 
	Let $(Y,\tau)$ be a cubic threefold with a non-Eckardt type involution and let $G:=\mathbb{Z}_2=\langle 1,\tau\rangle$. The fixed locus of $\tau$ consists of a pointwise fixed line $L$ and a plane cubic curve $C$: $Y^\tau=C\coprod L$. According to the seminal work \cite{BKR2001mckay}, there is an equivalence between the equivariant derived category $D^{b}_{\mathbb{Z}_2}(Y)$ and the derived category $D^{b}(Z)$, where $Z$ is an irreducible component of the $G\text{-Hilbert}$ scheme $G\text{-Hilb}_{\mathbb{C}}(Y)$, and is birational to $Y/\tau$. Following the strategy in \cite{Hu2023Equivariant}, we give an explicit description of the irreducible component $Z$. To be precise, we show that the irreducible component $Z \cong Bl_{\widetilde{C}}(\mathbb{P}^{2}\times\mathbb{P}^{1})$ (cf.~Proposition \ref{Theorem about the important component of G hilbert scheme}) where $\widetilde{C}$ coincides with the double cover of the cubic component $C\subset Y^\tau$ in \cite{casalaina2022moduli} (recall that $\widetilde{C}$ also parametrizes $\tau$-invariant lines in $Y$ which are not pointwise fixed). From this, we derive a semi-orthogonal decomposition of the equivariant category $D^{b}_{\mathbb{Z}_2}(Y)$ (cf.~Theorem \ref{thm_description_equivariant_derived_category})
	$$D^b_{\mathbb{Z}_2}(Y)\simeq D^b(Z)=\langle D^b(\widetilde{C}),E_{ij}\rangle$$ where $\widetilde{C}$ is described above and $\{E_{ij}\}$ with $1\leq i\leq 3,1\leq j\leq 2$ is an exceptional collection of line bundles on $\mathbb{P}^{2}\times\mathbb{P}^{1}$. Based on this semi-orthogonal decomposition, we then explore the Hodge theory of the equivariant Kuznetsov component $\Ku_{\mathbb{Z}_2}(Y)$. Specifically, we generalize the intermediate Jacobian construction in \cite{perry2020integral} to a smooth and proper $dg$ category (see Definition \ref{generalization_of_Alex_J}) and then adapt it to $\Ku_{\mathbb{Z}_2}(Y)$. This, together with the observation that the topological K-group $K^{top}_{1}(\Ku_{\mathbb{Z}_2}(Y))$ inherits a polarized Hodge structure of weight $1$ from $\mathrm{H}^1(\widetilde{C},\mathbb{Z})$, allows us to obtain an isomorphism $J(\Ku_{\mathbb{Z}_2}(Y))\cong J(\widetilde{C})$ of principally polarized abelian varieties (see Proposition \ref{Lemma_intermediate_Jacobian_Ku_same_curve}). By \cite{CDC05}, a general bielliptic curve of genus $4$ admits a unique bielliptic structure; note also that the invariant part $J(Y)$ is isomorphic to the dual abelian variety of $P(\widetilde{C}, C)$. A first proof of Theorem \ref{main_theorem_first_equivariant_categorical_Torelli} can then be completed (for general cubic threefolds $(Y,\tau)$ with a non-Eckardt type involution) using the results in \cite{casalaina2022moduli}.

\subsection{Bridgeland moduli spaces on the equivariant Kuznetsov components}\label{section_Bridgeland}
	Notation as in the previous subsection. Our second approach to proving Theorem \ref{main_theorem_first_equivariant_categorical_Torelli} is to use stability conditions on the equivariant Kuznetsov components.  Specifically, a Serre-invariant stability condition on $\Ku(Y)$ induces a (unique) stability condition $\sigma^{\mathbb{Z}_2}$ on $\Ku_{\mathbb{Z}_2}(Y)$. By studying the moduli space of $\sigma^{\mathbb{Z}_2}$-stable objects, we reconstruct in Proposition \ref{invariant_lines_as_Bridgeland_moduli_space} two copies of fixed locus $F(Y)^{\tau}$ of the Fano surface of lines in the non-Eckardt type cubic threefold $(Y,\tau)$. By a result in \cite{casalaina2022moduli}, $F(Y)^{\tau}$ consists of two components: a point corresponding to the pointwise fixed line $L$ and a curve parametrizing other $\tau$-invariant lines; note that the curve is isomorphic to the double cover $\widetilde{C}$ of the cubic curve $C\subset Y^\tau$. An  argument similar to the one described in the previous subsection then allows us to complete the proof.
	
\subsection{Chow theory of the equivariant Kuznetsov components}\label{section_Chow}
	The third proof of a slightly weaker version of Theorem \ref{main_theorem_first_equivariant_categorical_Torelli} we  complete is via the Chow theory of the equivariant Kuznetsov components. More precisely, we show that the Fourier--Mukai type equivalence $\Phi:\Ku_{\mathbb{Z}_2}(Y)\simeq\Ku_{\mathbb{Z}_2}(Y')$ in Theorem \ref{main_theorem_first_equivariant_categorical_Torelli} induces an isomorphism between the groups $A_{\mathbb{Z}_2,\mathbb{Q}}^2(Y)$ and $A_{\mathbb{Z}_2,\mathbb{Q}}^2(Y')$ of $\mathbb{Z}_2$-invariant algebraically trivial cycles with rational coefficients. Composing with the Abel--Jacobi map, one obtains an isogeny between the invariant parts $J(Y)^\tau$ and $J(Y')^{\tau'}$ of the intermediate Jacobians. It then follows from \cite{casalaina2022moduli} that $(Y,\tau)\cong (Y',\tau')$ when both of them are very general.
	
	 To obtain the above results, we invoke the techniques developed in \cite{kuznetsovnonclodedfield2021}. While we are only able to prove the equivariant categorical Torelli for very general non-Eckardt cubic threefolds via this approach, similar arguments allow us to prove the following theorem in a more general setting.
\begin{theorem}[=Theorem \ref{thm_general_result_chow_isom}]\label{main_theorem_chow}
	Let $X$ and $Y$ denote smooth projective rationally connected threefolds, each admitting an action of a finite group $G$.  Assume that both $D^b(X)$ and $D^b(Y)$ admit semi-orthogonal decompositions
	$$D^b(X)=\langle \mathcal A,{}^\perp\mathcal{A}\rangle,\,\,\,\,\,\,
D^b(Y)=\langle \mathcal{B},{}^{\perp}\mathcal{B}\rangle$$
such that ${{}^{\perp}\mathcal{A}}$ and ${}^{\perp}\mathcal{B}$ are generated by exceptional collections. 
 If there exists a Fourier--Mukai type equivalence $\Phi$ with kernel $P$ between the equivariant components
	$$\Phi:\mathcal A_G\stackrel{\simeq}{\to}\mathcal B_G,$$ 
then we have an isomorphism between the groups of $G$-invariant algebraically trivial cycles with coefficients in $\mathbb{Z}[1/2]$ (where $P'$ denotes the kernel of inverse Fourier-Mukai functor of $\Phi$)
	$$\xymatrix{A^{2}_{G, \mathbb{Z}[1/2]}(X)\ar@/^/[r]^{\frac{c_{3}(P)}{2}}& A^{2}_{G, \mathbb{Z}[1/2]}(Y)\ar@/^/[l]^{\frac{c_{3}(P')}{2}} 
  }.$$

	Consequently, there is an isogeny $J_G(X)\to J_G(Y)$ between the $G$-invariant part of intermediate Jacobians; moreover, the order of the kernel is divisible at most by a power of $2$.
 %the primes dividing $|G|$ and $2$.
\end{theorem}

\begin{remark}\label{remark_alternative_approach}
	It is suggested by Sasha Kuznetsov and Xiaolei Zhao that if one further assumes that the equivalence $\Phi:\Ku_{\mathbb{Z}_2}(Y)\simeq\Ku_{\mathbb{Z}_2}(Y')$ in Theorem \ref{main_theorem_first_equivariant_categorical_Torelli} preserves the residual actions of the dual group $\widehat{\mathbb{Z}_2}$ (which is given by tensoring with the characters of $\mathbb{Z}_2$), then the genericity assumption in Theorem~\ref{main_theorem_first_equivariant_categorical_Torelli} can be removed. We provide a proof in the appendix (Theorem~\ref{thm_equivariant_cat_Torelli}), which is more category theoretical and relies on Elagin's reconstruction theorem \cite[Theorem 4.2]{elagin2014equivariant} and the classification of auto-equivalences of Kuznetsov components of cubic threefolds in \cite{feyzbakhsh2023new} and \cite{liu2024fourier}. 
    
    %Indeed, since the equivalence $\Phi$ commutes with the residual actions of $\widehat{\mathbb{Z}_2}$, it induces an equivalence $\Psi:\Ku_{\mathbb{Z}_2}(Y)^{\widehat{\mathbb{Z}_2}}\simeq\Ku_{\mathbb{Z}_2}(Y')^{\widehat{\mathbb{Z}_2}}$ of the invariant categories, which is also compatible with the actions of $\widehat{\widehat{\mathbb{Z}_2}}\cong\mathbb{Z}_2$ (cf.~the proof of \cite[Lemma 4.9]{bayer2023kuznetsov}). By Elagin's reconstruction theorem, it also holds that $\Ku_{\mathbb{Z}_2}(Y)^{\widehat{\mathbb{Z}_2}}\simeq\Ku(Y)$. Thus we get an equivalence $\Psi:\Ku(Y)\simeq\Ku(Y')$ satisfying $\Psi\circ\tau^*\cong\tau'^*\circ\Psi$, where $\tau^*$ and $\tau'^*$ are the auto-equivalences of $\Ku(Y)$ and $\Ku(Y')$ coming from the geometric involutions of $Y$ and $Y'$, respectively. As a result of \cite[Corollary 6.11]{Li2024Higher}, $\Psi$ is of Fourier-Mukai type. Then by \cite[Theorem 4.23]{liu2024fourier} or \cite[Corollary 8.4]{feyzbakhsh2023new}, there exists a unique isomorphism $f:Y\cong Y'$ such that $f_*\circ\tau_*\cong\tau'_*\circ f_*$, and hence we get $(f\circ\tau\circ f^{-1})_*\cong\tau'_*$, where both sides are auto-equivalences of $\Ku(Y')$. Using \cite[Remark 8.5]{feyzbakhsh2023new}, we obtain that $f\circ\tau\cong\tau'\circ f$ which, together with the categorical Torelli theorem for cubic threefolds, allows us to conclude that $(Y,\tau)\cong (Y',\tau')$. 
\end{remark}

\subsection{Infinitesimal categorical Torelli theorem for equivariant Kuznetsov components}
	In the remaining part of the paper, we focus on the infinitesimal categorical Torelli problem. Let us discuss some of the key ingredients in the proof, which we hope will be of independent interest. Let $X$ be a smooth projective variety. By comparing the algebra structures on the Hochschild cohomology (respectively, the module structures of the Hochschild homology over the Hochschild cohomology) for $D^b(X)$ and for a semi-orthogonal component $\cA_1$ of $D^b(X)$, we establish the following compatibility theorem.     
	
\begin{theorem}[=Theorem \ref{compatibilitymain}]\label{theorem_second_main_result}
Let $\Phi_H$ (respectively, $\Phi_{1H}$) be an auto-equivalence of $D^b(X)$ (respectively, of a semi-orthogonal component $\cA_1$), and let $\Phi_{H\ast}$ and $\Phi_{1H\ast}$ be the induced automorphism on Hochschild cohomology $\mathrm{HH}^{\ast}(X)$ and $\mathrm{HH}^{\ast}(\cA_1)$ respectively.  Denote by $j_1^*$ the projection functor from $D^b(X\times X)$ to its semi-orthogonal component $\cA_1\boxtimes\mathcal{B}^{\vee}_1$. Then $\Phi_{H\ast}$ and $\Phi_{1H\ast}$ are compatible with both the algebra structures of the Hochschild cohomology and the module structures of the Hochschild homology over the Hochschild cohomology. In other words, the following diagrams are commutative.
\[\begin{tikzcd}
	&& \mathrm{HH}^{\ast}(\cA_1)\times\mathrm{HH}^{\ast}(\cA_1) & \mathrm{HH}^{\ast}(\cA_1)\times\mathrm{HH}^{\ast}(\cA_1) \\
	\mathrm{HH}^{\ast}(X)\times\mathrm{HH}^{\ast}(X) & \mathrm{HH}^{\ast}(X)\times\mathrm{HH}^{\ast}(X) & \mathrm{HH}^{\ast}(\cA_1) & \mathrm{HH}^{\ast}(\cA_1) \\
	\mathrm{HH}^{\ast}(X) & \mathrm{HH}^{\ast}(X)
	\arrow[from=2-1, to=1-3, "j_1^*"]
	\arrow[from=2-2, to=1-4, "j_1^*"]
	\arrow[dashed, from=3-1, to=2-3, "j_1^*"]
	\arrow[from=3-2, to=2-4, "j_1^*"]
	\arrow[from=2-1, to=3-1, "\cup"]
	\arrow[from=2-2, to=3-2, "\cup"]
	\arrow[from=2-1, to=2-2, "\Phi_{H_*}"]
	\arrow[from=3-1, to=3-2, "\Phi_{H_*}"]
	\arrow[from=1-3, to=1-4, "\Phi_{1H_*}"]
	\arrow[from=2-3, to=2-4, "\Phi_{1H_*}"]
	\arrow[from=1-4, to=2-4, "\cup"]
	\arrow[dashed, from=1-3, to=2-3, "\cup"]
\end{tikzcd}\]

\[\begin{tikzcd}
	&& \mathrm{HH}^{\ast}(\cA_1)\times\mathrm{HH}_*(\cA_1) & \mathrm{HH}^{\ast}(\cA_1)\times\mathrm{HH}_*(\cA_1) \\
	\mathrm{HH}^{\ast}(X)\times\mathrm{HH}_*(X) & \mathrm{HH}^{\ast}(X)\times\mathrm{HH}_*(X) & \mathrm{HH}_*(\cA_1) & \mathrm{HH}_*(\cA_1) \\
	\mathrm{HH}_*(X) & \mathrm{HH}_*(X)
	\arrow[from=2-1, to=1-3, "j_1^*"]
	\arrow[from=2-2, to=1-4, "j_1^*"]
	\arrow[dashed, from=3-1, to=2-3, "j_1^*"]
	\arrow[from=3-2, to=2-4, "j_1^*"]
	\arrow[from=2-1, to=3-1, "\cap"]
	\arrow[from=2-2, to=3-2, "\cap"]
	\arrow[from=2-1, to=2-2, "\Phi_{H_*}"]
	\arrow[from=3-1, to=3-2, "\Phi_{H_*}"]
	\arrow[from=1-3, to=1-4, "\Phi_{1H_*}"]
	\arrow[from=2-3, to=2-4, "\Phi_{1H_*}"]
	\arrow[from=1-4, to=2-4, "\cap"]
	\arrow[dashed, from=1-3, to=2-3, "\cap"]
\end{tikzcd}\] 
 \end{theorem}
 
	We also investigate the compatibility problem between the twisted Hochschild-Kostant-Rosenberg ($\operatorname{HKR}$) isomorphism (cf.~\cite{cualduararu2005mukai} and \cite{calaque2012cualduararu}) and the actions on Hochschild cohomology and on singular cohomology. Namely, a result in \cite{macri2009infinitesimal} implies that the $\operatorname{HKR}$ isomorphism is compatible with the actions induced by a Fourier-Mukai transform for Hochschild homology and singular cohomology. We extend this result to Hochschild cohomology under the assumption that the action is induced by a geometric automorphism (cf.~Theorem \ref{IKequivariant}).
	
	The above results allow us to prove the following theorem which is an equivariant version of the results  (especially Theorems 1.1 and 1.4) in \cite{jacovskis2023infinitesimal}. Noting that the injectivity of $dp$ implies the injectivity of $\eta$, we regard Theorem \ref{main_theorem_comm_diag_equivariant} as (an equivariant version of) the infinitesimal categorical Torelli theorem. Indeed, let $X$ be a smooth projective variety. Although it is difficult to make sense of  ``categorical period maps", one could still define a map $H^1(X,T_X)\rightarrow \HH^2(\mathcal{K}u(X))$ as in \cite{jacovskis2023infinitesimal}. Since first order deformations of $\Ku(X)\subset D^b(X)$ are parametrized by the second Hochschild cohomology $\mathrm{HH}^2(\Ku(X))$, the map $H^1(X,T_X)\rightarrow \HH^2(\mathcal{K}u(X))$ can be interpreted as the infinitesimal categorical period map. Taking the invariant parts, we think of the injectivity of the map $\eta$ in Theorem \ref{main_theorem_comm_diag_equivariant} as an equivariant infinitesimal categorical Torelli theorem.

\begin{theorem}[=Theorem \ref{main_theorem_compatibility}, Corollary \ref{cor_main_theorem_compatibility}]\label{main_theorem_comm_diag_equivariant}
	Let $X$ be a smooth projective variety with a biregular involution $\tau$. Suppose that $D^b(X)$ admits a semi-orthogonal decomposition 
	$$D^b(X)=\langle\Ku(X),E_1,\ldots,E_n\rangle.$$ 
Also assume that the involution $\tau$ induces a $\mathbb{Z}_2$-action $\tau^*$ on $\Ku(X)$ which preserves each semi-orthogonal component. Then we have a commutative diagram as follows. 
    $$\xymatrix@C=2.5cm{\mathrm{HH}^{2}(\Ku(X))^\tau\ar[r]^(0.4){\gamma}&Hom(H\Omega_{-1}(X)^{\tau},H\Omega_{1}(X)^{\tau})\\
H^{1}(X,T_{X})^{\tau}\ar[u]^{\eta}\ar[ur]^(0.4){dp}& }$$

	Moreover, if $(Y,\tau)$ is a general non-Eckardt cubic threefold, then  the map $H^{1}(Y,T_{Y})^{\tau}\stackrel{\eta}{\to} \mathrm{HH}^{2}(\Ku(Y))^{\tau}\subset \HH^{2}(\Ku_{\mathbb{Z}_2}(Y))$ is injective.
\end{theorem}

\subsection{Organization of the article}
	In Section \ref{section_equivariant_SOD}, we recall some background material on semi-orthogonal decompositions and equivariant categories. In Section \ref{section_cubic_involution}, we review results on cubic threefolds with a non-Eckartd type involution. In Section \ref{section_description_equivariant_Ku}, we give an explicit description of the equivariant derived categories for non-Eckardt cubic threefolds. In Section \ref{section_equivariant_CT}, we explore the noncommutative Hodge theory of the equivariant Kuznetsov components and give the first proof of Theorem \ref{main_theorem_first_equivariant_categorical_Torelli}. In Section \ref{section_reconstruct_Bridgeland_moduli_space}, we investigate the equivariant Kuzentsov components via Bridgeland moduli which allows us to give a  second proof of Theorem \ref{main_theorem_first_equivariant_categorical_Torelli}. In Section \ref{section_chow_theory_Kuznetsov_components}, we consider the Chow theoretic perspective, and show that the Fourier--Mukai equivalence in Theorem \ref{main_theorem_first_equivariant_categorical_Torelli} induces an isomorphism between the groups of invariant algebraically trivial cycles. As a result, we provide the third proof (but only for very general non-Eckardt cubic threefolds). In Section \ref{section_inf_cat_Torelli}, we discuss the proofs of Theorem \ref{theorem_second_main_result} and Theorem \ref{main_theorem_comm_diag_equivariant} regarding infinitesimal categorical Torelli theorems for the equivariant Kuznetsov components. 
	
	Throughout the article, we work over the field of complex numbers $\mathbb{C}$. The main reason is that we need to use the results in \cite{casalaina2022moduli} which are only stated over $\mathbb{C}$; note however that most results in Sections \ref{section_equivariant_SOD}, \ref{section_equivariant_CT}, \ref{section_reconstruct_Bridgeland_moduli_space}, \ref{section_chow_theory_Kuznetsov_components}, and \ref{section_inf_cat_Torelli} hold for any algebraically closed field of characteristic $0$. Furthermore, the  results in \S \ref{section_chow_theory_Kuznetsov_components}, with the exception of those relying on  \cite{casalaina2022moduli}, can be formulated over any field, after replacing the intermediate Jacobian with the algebraic representative, but we do not pursue this here. 

\subsection{Acknowledgement}
	We would like to thank Nicolas Addington, Alexey Elagin, Daniel Huybrechts, Shengxuan Liu, Lisa Marquand, Alexander Perry, Laura Pertusi, Paolo Stellari and Junwu Tu for useful discussions. We are indebted to Sasha Kuznetsov and Xiaolei Zhao for pointing out Remark \ref{remark_alternative_approach} to us and for related discussions. Special thanks also go to Sasha Kuznetsov for multiple valuable suggestions. We are grateful to Augustinas Jacovskis as well for helping us in drawing several commutative diagrams. S.C.-M. is supported in part by a grant from the Simons Foundation (581058). X.H. was funded by the Deutsche Forschungsgemeinschaft (DFG, German Research Foundation) -- 516701553. S.Z. is supported by ANR project FanoHK, grant ANR-20-CE40-0023, Deutsche Forschungsgemeinschaft under Germany's Excellence Strategy-EXC-2047$/$1-390685813 and NSF grant No. DMS-1928930 (while he is working at the Simons Laufer Mathematical Sciences Institute in Berkeley, California). S.Z. is also partially supported by GSSCU 2021092. Z.Z. is supported in part by NSFC grant 12201406. Part of the work was completed when some of the authors visited Yau Mathematical Science Center at Tsinghua University, Max-Planck Institute for Mathematics, and Hausdorff Research Institute for Mathematics. We are grateful for their excellent hospitality and support.

\section{Semi-orthogonal decompositions and equivariant derived categories}\label{section_equivariant_SOD}
	In this section, let us briefly review some background material on semi-orthogonal decompositions and equivariant derived categories. We refer the reader to \cite{bondal1989representable}, \cite{bondal1989representations} and \cite{beckmann2023equivariant} for more details. We also specify the geometric situation when there is a smooth cubic threefold $Y$ admitting an involution $\tau$ and show that its equivariant derived category $D^b_{\mathbb{Z}_2}(Y)$ admits a semi-orthogonal decomposition with one triangulated subcategory being the equivariant Kuznetsov component $\Ku_{\mathbb{Z}_2}(Y)$. A more concrete description of the semi-orthogonal decomposition of  $D^b_{\mathbb{Z}_2}(Y)$ and of the equivariant Kuznetsov component $\Ku_{\mathbb{Z}_2}(Y)$ will be given in the later sections.

\subsection{Semi-orthogonal decompositions}
	Let $\mathcal{D}$ be a $\mathbb{C}$-linear triangulated category. A \emph{semi-orthogonal decomposition} of $\mathcal{D}$ consists of a collection $\mathcal{C}_1,\ldots,\mathcal{C}_m$ of full triangulated subcategories satisfying that 
\begin{enumerate}
    \item $\mathrm{Hom}(F,G)=0$ for every $F\in\mathcal{C}_i$, $G\in\mathcal{C}_j$ and $i>j$;
    \item for every $F\in\mathcal{D}$, there is a sequence of morphisms 
    $$0=F_m\rightarrow F_{m-1}\rightarrow\ldots\rightarrow F_1\rightarrow F_0=F$$
    such that the cone of $F_i\rightarrow F_{i-1}$ is in $\mathcal{C}_i$ for $1\leq i\leq m$.
\end{enumerate}
As usual, we use $\mathcal{D}=\langle\mathcal{C}_1,\ldots,\mathcal{C}_m\rangle$ to denote a semi-orthogonal decomposition. Recall that an object $E\in\mathcal{D}$ is called \emph{exceptional} if $\mathrm{RHom}(E,E)=\mathbb{C}[0]$. Given an \emph{exceptional collection} $E_1,\ldots,E_m$ of $\mathcal{D}$ (that is, a sequence of exceptional objects $E_1,\ldots,E_m\subset\mathcal{D}$ satisfying $\mathrm{RHom}(E_i,E_j)=0$ for $i>j$), we get the following semi-orthogonal decomposition 
$$\mathcal{D}=\langle\mathcal{C},E_1,\ldots,E_m\rangle$$ where $\mathcal{C}:=\langle E_1,\ldots,E_m\rangle^{\perp}=\{F\in\mathcal{D}\mid\mathrm{RHom}(E_i,F)=0$ for all $1\leq i\leq m\}$ denotes the right orthogonal subcategory of the exceptional collection. 

\subsection{Categorical actions and equivariant categories}
	Let us now recall some basic definitions of categorical actions and equivariant categories following \cite{beckmann2023equivariant}. 

\begin{definition}{(\cite[Definition 2.1]{beckmann2023equivariant})}\label{def_group_action_category}
	Let $G$ be a finite group and let $\mathcal{D}$ be a category. An \emph{action} $(\rho, \theta)$ of $G$ on $\mathcal{D}$ consists of
\begin{itemize}
\item for every $g \in G$ an auto-equivalence $\rho_g \colon \mathcal{D} \to \mathcal{D}$ (we will sometimes write $g$ for $\rho_g$),
\item for every pair $g,h \in G$ an isomorphism of functors $\theta_{g,h} \colon \rho_{g} \circ \rho_h \to \rho_{gh}$
\end{itemize}
such that for all triples $g,h,k \in G$ we have the commutative diagram.
\begin{equation} \label{associativity}
\begin{tikzcd}[row sep=large, column sep = large]
\rho_{g} \rho_{h} \rho_{k} \ar{r}{\rho_g  \theta_{h,k}} \ar{d}{\theta_{g,h} \rho_k} & \rho_{g} \rho_{hk} \ar{d}{\theta_{g,hk}} \\
\rho_{gh} \rho_{k} \ar{r}{\theta_{gh,k}} & \rho_{ghk}
\end{tikzcd}
\end{equation}
\end{definition}

\begin{definition}{(\cite[Definition 3.1]{beckmann2023equivariant})}\label{def_equivariant_category}
	Let $(\rho,\theta)$ be an action of a finite group $G$ on a  $\mathbb{C}$-linear category $\mathcal{D}$. The \emph{equivariant category} $\mathcal{D}_G$ is defined as follows:
\begin{itemize}
\item Objects of $\mathcal{D}_G$ are pairs $(E,\phi)$ where $E$ is an object in $\mathcal{D}$ and $\phi = (\phi_g \colon E \to \rho_{g} E)_{g \in G}$ is a family of isomorphisms
such that
\begin{equation} \label{compatibility}
\begin{tikzcd}
E \ar[bend right]{rrr}{\phi_{gh}} \ar{r}{\phi_g} & \rho_{g}E \ar{r}{\rho_g \phi_h} &  \rho_{g} \rho_{h} E \ar{r}{\theta_{g,h}^E} & \rho_{gh} E
\end{tikzcd}
\end{equation} 
commutes for all $g,h \in G$.
\item A morphism from $(E,\phi)$ to $(E', \phi')$ is a morphism $f \colon E \to E'$ in $\mathcal{D}$
which commutes with linearizations, i.e.\ such that
\[
\begin{tikzcd}
E \ar{r}{f} \ar{d}{\phi_g} & E' \ar{d}{\phi'_g} \\
g E \ar{r}{\rho_gf} & g E'
\end{tikzcd}
\]
commutes for every $g \in G$.
\end{itemize}
\end{definition}

\begin{remark}
	As mentioned in \cite[p.36]{beckmann2023equivariant}, for any objects $(E,\phi)$ and $(E',\phi')$ in $\mathcal{D}_G$ there exists an action of $G$ on $\Hom_{\mathcal{D}}(E,E')$ via 
\[ f \mapsto (\phi'_g)^{-1} \circ \rho_g(f) \circ \phi_g. \]
Then it holds that 
\[ \Hom_{\mathcal{D}_G}( (E,\phi), (E',\phi') ) = \Hom_{\mathcal{D}}(E,E')^G. \]
\end{remark}

	Now we focus on the following geometric situation. Let $X$ be a smooth projective variety and let $G$ be a finite group acting on $X$. Suppose that $D^b(X)$ admits a semi-orthogonal decomposition 
	$$D^b(X)=\langle\mathcal{A}_1,\mathcal{A}_2,\ldots,\mathcal{A}_m\rangle.$$ 
Also assume that $D^b(X)$ and each $\mathcal{A}_i$ admits an action of $G$ (in the sense of Definition \ref{def_group_action_category}). Then by \cite[Theorem 6.3]{elagin2011cohomological}, we get a semi-orthogonal decomposition 
	$$D^b_G(X)=\langle\mathcal{A}_{1G}, \mathcal{A}_{2G},\ldots,\mathcal{A}_{mG}\rangle$$
of the equivariant derived category $D^b_G(X)$. 

	 More specifically, let $Y$ be a smooth cubic threefold with an involution $\tau$. Recall that the derived category $D^b(Y)$ admits a semi-orthogonal decomposition $$D^b(Y)=\langle\Ku(Y),\oh_Y,\oh_Y(1)\rangle$$ where $\Ku(Y)$ denotes the right orthogonal subcategory of the exceptional collection $\oh_Y$ and $\oh_Y(1)$, and is called the \emph{Kuznetsov component} of $Y$. Now let us introduce the equivariant derived category $D^b_{\mathbb{Z}_2}(Y)$ and the \emph{equivariant Kuznetsov component} $\Ku_{\mathbb{Z}_2}(Y)$ for the action by the group $\mathbb{Z}_2=\langle 1,\tau\rangle$.

\begin{lemma}\label{cubic_threefold_equivariant_SOD}
	Let $(Y,\tau)$ be as above. There exists a semi-orthogonal decomposition 
$$D^b_{\mathbb{Z}_2}(Y)=\langle\Ku_{\mathbb{Z}_2}(Y), \oh_Y, \oh_Y\otimes\chi_1, \oh_Y(1),\oh_Y(1)\otimes\chi_1\rangle$$ where $\chi_1$ denotes the sign character of the group $\mathbb{Z}_2=\langle 1,\tau\rangle$. 
\end{lemma}

\begin{proof}
 	Consider the semi-orthogonal decomposition of $Y$: $D^b(Y)=\langle\Ku(Y),\oh_Y,\oh_Y(1)\rangle$. Since $\tau$ is a geometric involution of $Y$, it clearly preserves $\Ku(Y)$, $\oh_Y$ and $\oh_Y(1)$. It is straightforward to verify that both $D^b(Y)$ and $\Ku(Y)$ admit an action by $\mathbb{Z}_2=\langle 1,\tau\rangle$ as in Definition \ref{def_group_action_category} (see also Remark \ref{remk_existence_action}).  By \cite[Theorem 1.6]{elagin2011cohomological}, we get the following semi-orthogonal decomposition of the equivariant derived category $D^b_{\mathbb{Z}_2}(Y)$:
	$$D^b_{\mathbb{Z}_2}(Y)=\langle\Ku_{\mathbb{Z}_2}(Y),\langle\oh_Y\rangle^{\tau},\langle\oh_Y(1)\rangle^{\tau}\rangle.$$
The result then follows from \cite[Proposition 3.3]{KP2017}.
\end{proof}

%\begin{remark} 
%When the action $\tau$ is of non-Eckardt type, one could also use \cite[Theorem 1.1]{li2022derived} and \cite[Corollary 4.10]{beckmann2023equivariant} to conclude there is an induced action of $\mathbb{Z}_2=\langle 1,\tau\rangle$ on the Kuznetsov component $\Ku(Y)$. See Remark \ref{remk_existence_action}.}
%\end{remark}

\section{Cubic threefolds with an involution of non-Eckardt type}\label{section_cubic_involution}
	In this section, we recall some related results on cubic threefolds with a non-Eckardt type involution following \cite{casalaina2022moduli}.

	Let $Y$ be a smooth cubic threefold with an involution $\tau$. Then the fixed locus of $\tau$ is either a point and a cubic surface (in which case we say $\tau$ is \emph{of Eckardt type}) or a line disjoint union with a plane cubic curve (in which case we say $\tau$ is \emph{of non-Eckardt type}). These cubic threefolds have been studied in detail in \cite{casalaina2021moduli} and \cite{casalaina2022moduli} respectively.

	In this paper, we will focus on cubic threefolds $Y$ with a non-Eckardt type involution $\tau$. After a linear change of coordinates, the equation of $Y$ is given by 
\begin{equation}\label{eqn_nEkcardt_Y}
	x_0q_0(x_3,x_4)+x_1q_1(x_3,x_4)+x_2q_2(x_3,x_4)+g(x_0,x_1,x_2)=0
\end{equation}
where each polynomial $q_i(x_3,x_4)$ is homogeneous of degree $2$ and $g(x_0,x_1,x_2)$ is a homogeneous cubic polynomial. The non-Eckardt type involution
\begin{equation}\label{eqn_nEcardt_involution}
    \tau:[x_0,x_1,x_2,x_3,x_4]\mapsto [x_0,x_1,x_2,-x_3,-x_4]
\end{equation}
fixes pointwise the line 
    $$L:=V(x_0,x_1,x_2)\subset Y$$
and the plane cubic curve (which is smooth)
	$$C:=V(g(x_0,x_1,x_2),x_3,x_4)\subset Y.$$

	Let $J(Y)$ be the intermediate Jacobian of $Y$. By abuse of notation, we also use $\tau$ to denote the involution on $J(Y)$ induced by the non-Eckardt involution on $Y$. Define the invariant part $J(Y)^\tau$ by 
	$$J(Y)^\tau:=\im(1+\tau)$$
which admits a polarization of type $(1,2,2)$ from the principal polarization on $J(Y)$. By \cite[Theorem 3.1]{casalaina2022moduli}, we have the following global Torelli theorem.

\begin{theorem}{(\cite[Theorem 3.1]{casalaina2022moduli})}\label{thm_global_torelli_nEckardt}
	Let $(Y,\tau)$ and $(Y',\tau')$ be cubic threefolds with an involution of non-Eckardt type. Suppose that $J(Y)^\tau\cong J(Y')^{\tau'}$ as polarized abelian varieties, then $(Y,\tau)\cong (Y',\tau')$. 
\end{theorem}

	The above global Torelli theorem can also be rephrased as follows. Let $\mathcal{M}$ be the moduli space of cubic threefolds with an involution of non-Eckardt type. Denote the moduli space of abelian threefolds with a polarization of type $(1,2,2)$ by $\mathcal{A}_3^{(1,1,2)}$. Define the period map as follows:
	$$\mathcal{P}:\mathcal{M}\to \mathcal{A}_3^{(1,2,2)},\,\,\,\,\,\,Y\to J(Y)^\tau.$$
Then the above period map $\mathcal{P}$ is injective.

	To prove the above the theorem, we project $Y$ from the $\tau$-fixed line $L\subset Y$ to the complementary plane $\mathbb{P}^{2}_{[x_{0}:x_{1}:x_{2}]}$. The discriminant curve consists of two irreducible components meeting transversely: the cubic curve $C$ and a conic curve $Q$. The discriminant double cover is given by 
	$$\widetilde{C}\coprod \widetilde{Q}\to C\coprod Q$$ 
where $\widetilde{C}$ (respectively, $\widetilde{Q}$) is a double cover of the cubic $C$ (respectively, the conic $Q$) branched at the six intersection	points $C\cap Q$ (cf.~\cite[Lemma 2.1, Proposition 2.2]{casalaina2022moduli}). As discussed in \cite[\S1.3]{casalaina2022moduli}, the bielliptic curve $\widetilde{C}$ is smooth of genus $4$ and parameterizes $\tau$-invariant lines in $Y$ which are not pointwise fixed. By \cite[Theorem 2.5]{casalaina2022moduli}, the invariant part $J(Y)^\tau$ is isomorphic to the dual abelian variety of the Prym variety $P(\widetilde{C}, C)$ associated with the discriminant double restricted to the cubic component $\widetilde{C}\to C$:
	$$J(Y)^\tau\cong (P(\widetilde{C},C))^\vee.$$
Another important ingredient for the proof is the global Torelli theorem for the Prym map $\mathcal{P}_{1,6}:\mathcal{R}_{1,6}\to \mathcal{A}_3^{(1,1,2)}$ proved by Ikeda \cite[Theorem 1.2]{MR4156417} and by Naranjo and Ortega \cite[Theorem 1.1]{MR4484544}.

	The infinitesimal Torelli theorem also holds for the period map $\mathcal{P}:\mathcal{M}\to \mathcal{A}_3^{(1,2,2)}$ over an open subset $\mathcal{M}_0\subset \mathcal{M}$ (cf.~\cite[Proposition 3.4]{casalaina2022moduli}). Specifically, let $\mathcal{M}_0$ denote the complement of the locus parameterizing cubic threefolds with a non-Eckardt involution whose equation can be written in the following form
	$$\ell_1(x_0,x_1,x_2)x_3^2+\ell_2(x_0,x_1,x_2)x_4^2+g(x_0,x_1,x_2)=0.$$
Then we have the following result.

\begin{proposition}{(\cite[Proposition 3.4]{casalaina2022moduli})}\label{prop_inf_torelli_nEckardt}
    The differential $d\mathcal{P}$ is an isomorphism at every point of $\mathcal{M}_0\subset \mathcal{M}$. In particular, $\mathcal{P}|_{\mathcal{M}_0}:\mathcal{M}_0\to \mathcal{A}_3^{(1,2,2)}$ is an embedding.
\end{proposition}

\section{Semi-orthogonal decompositions of the equivariant derived categories}\label{section_description_equivariant_Ku}
	Let $Y$ be a cubic threefold with a non-Eckardt type involution $\tau$ as in Section \ref{section_cubic_involution}. We have defined the equivariant derived category $D^b_{\mathbb{Z}_2}(Y)$ for the action of $G:=\mathbb{Z}_2=\langle 1,\tau\rangle$ in Section \ref{section_equivariant_SOD}. In this section, we give an explicit description of $D^b_{\mathbb{Z}_2}(Y)$ via $G$-Hilbert schemes. More precisely, we will prove the following theorem.

\begin{theorem}\label{thm_description_equivariant_derived_category}
	The equivariant derived category $D^b_{\mathbb{Z}_2}(Y)$ admits a semi-orthogonal decomposition 
    $$D^b_{\mathbb{Z}_2}(Y)=\langle D^b(\widetilde{C}),E_{ij}\rangle,$$
where $E_{ij}= \mathcal{O}_{\mathbb{P}^{2}\times \mathbb{P}^{1}}(i,j)$ with $1\leq i\leq 3$ and $1\leq j\leq 2$ forms an exceptional collection of line bundles originating from $\mathbb{P}^2\times\mathbb{P}^1$, and $\widetilde{C}$ is the  double cover of the plane cubic curve $C$ in the fixed locus $Y^{\tau}$. 
\end{theorem}

	In this section, we always let $G=\mathbb{Z}_2=\langle 1,\tau\rangle$. To prove Theorem \ref{thm_description_equivariant_derived_category}, we shall employ the same strategy as in \cite{Hu2023Equivariant}. The idea is described as follows. Recall that a \emph{$G\text{-Hilbert}$ scheme} $G\text{-Hilb}_{\mathbb{C}}(Y)$ is a Hilbert scheme of $G$-clusters (see for instance \cite[\S13.1]{MR2244106} and the references therein). Let $Z$ be the irreducible component of $G\text{-Hilb}_{\mathbb{C}}(Y)$ containing the open subset of all reduced $G$-clusters; note that $Z$ admits a crepant resolution to $Y/G$. According to \cite[Theorem 1.1]{BKR2001mckay}, there is an equivalence between the equivariant derived category $D^{b}_{G}(Y)$ and the derived category $D^{b}(Z)$. By computing $Z$ explicitly, we derive a semi-orthogonal decomposition of the derived category $D^{b}(Z)$ with one component equivalent to the derived category of a certain curve. Upon further identification, we show that this curve is isomorphic to the double cover $\widetilde{C}$ of the plane cubic curve $C\subset Y^\tau$ described in Section \ref{section_cubic_involution}.
 
	In order to describe the irreducible component $Z\subset G\text{-Hilb}_{\mathbb{C}}(Y)$, let us study the $G-$Hilbert schemes $G\text{-Hilb}_{\mathbb{C}}(\mathbb{P}^{4})$ and $G\text{-Hilb}_{\mathbb{C}}(Y)$ explicitly. Consider the following involution on $\mathbb{P}^{4}$ whose restriction to $Y$ corresponds to the involution $\tau$ of non-Eckardt type 
	$$[x_{0}:x_{1}:x_{2}:x_{3}:x_{4}]\longrightarrow [x_{0}:x_{1}:x_{2}:-x_{3}:-x_{4}].$$
The fix locus of this involution in $\mathbb{P}^{4}$ is $\mathbb{P}^{2}\coprod \mathbb{P}^{1}$ where the coordinates of the fixed components are $[x_{0}:x_{1}:x_{2}:0:0]$ for $\mathbb{P}^{2}=\mathbb{P}^{2}_{[x_{0}:x_{1}:x_{2}]}$ (which intersects $Y$ along the cubic curve $C\subset Y^\tau$) and $[0:0:0:x_{3}:x_{4}]$ for $\mathbb{P}^{1}=\mathbb{P}^{1}_{[x_{3}:x_{4}]}$ (which is also the pointwise fixed line $L\subset Y^\tau$) respectively. Let $\mathbb{A}^{4}$ denote the open affine subset $(x_{0}\neq 0)$ of $\mathbb{P}^{4}$; correspondingly, one gets an open subset $G\text{-Hilb}_{\mathbb{C}}(\mathbb{A}^{4})$ of $G\text{-Hilb}_{\mathbb{C}}(\mathbb{P}^{4})$.  

\begin{lemma}\label{Lemma for G-Hilbert of A4}
	It holds that $G\text{-Hilb}_{\mathbb{C}}(\mathbb{A}^{4})\cong \mathbb{A}^{2}\times G\text{-Hilb}_{\mathbb{C}}(\mathbb{A}^{2})$ where for the first factor $\mathbb{A}^{2}=V(x_3, x_4)$ and for the second factor $\mathbb{A}^{2}=V(x_1, x_2)$ which is endowed with an induced involution given by $(x_{3}, x_{4})\longrightarrow (-x_{3}, -x_{4})$.
\end{lemma}
\begin{proof}
It follows directly from \cite[Proposition 1.4.4]{terouanne:tel-00006683} or \cite[Corrollary 4.24]{blume2007mckay}.
\end{proof}

We will also need the following result due to Blume.
\begin{proposition}{(\cite[Proposition 2.40]{blume2007mckay})}\label{Proposition by Blume}
	Let $G=\text{Spec }\mathbb{C}[x]/(x^{r}-1)$ be the group scheme of the $r$-th roots of unity, and let $V$ be an $n$-dimensional representation of $G$ over $\mathbb{C}$. Define $\pi: \mathbb{A}_{\mathbb{C}}(V)\longrightarrow X:= \mathbb{A}_{\mathbb{C}}(V)/G$ where $\mathbb{A}_{\mathbb{C}}(V)$ denotes the affine space $\mathbb{A}_{\mathbb{C}}^{n}$. Also denote the origin of $\mathbb{A}_{\mathbb{C}}(V)$ by $0$. Then there exists an isomorphism (over $\mathbb{C}$)
	$$G\text{-Hilb}_{\mathbb{C}}(\mathbb{A}_{\mathbb{C}}(V))\cong \text{Bl}_{0}X/G.$$
\end{proposition}

	Combining Lemma \ref{Lemma for G-Hilbert of A4} and Proposition \ref{Proposition by Blume}, we obtain the following description of $G\text{-Hilb}_{\mathbb{C}}(\mathbb{P}^{4})$.
\begin{lemma}\label{G-Hilbert scheme of P4}
The $G$-Hilbert scheme $G\text{-Hilb}_{\mathbb{C}}(\mathbb{P}^{4})$ is isomorphic to $\text{Bl}_{\mathbb{P}^{2}\coprod\mathbb{P}^{1}}(\mathbb{P}^{4})/G$.
\end{lemma}
\begin{proof}
According to \cite[Remark 4.19]{blume2007mckay}, if $X^{\prime}$ is a scheme endowed with a $G$-action over $S$ and $S^{\prime}$ is an $S$-scheme, then there is an isomorphism of $S^{\prime}$-functors
	$$(\underline{G\text{-Hilb}}_{S}(X^{\prime}))_{S^{\prime}}\cong \underline{G\text{-Hilb}}_{S^{\prime}}(X^{\prime}_{S^{\prime}}).$$
In addition, from \cite[Remark 4.22 (2)]{blume2007mckay} we deduce that 
\begin{center}
$G\text{-Hilb}_{\mathbb{C}}(\mathbb{P}^{4})\cong G\text{-Hilb}_{\mathbb{P}^{4}/G}(\mathbb{P}^{4})$ and $G\text{-Hilb}_{\mathbb{C}}(\mathbb{A}^{4})\cong G\text{-Hilb}_{\mathbb{A}^{4}/G}(\mathbb{A}^{4})$.
\end{center}
Thus, we have the following Cartesian diagram:
$$\begin{tikzcd}
G\text{-Hilb}_{\mathbb{C}}(\mathbb{A}^{4})\cong G\text{-Hilb}_{\mathbb{A}^{4}/G}(\mathbb{A}^{4}) \arrow[d, ""] \arrow[r, "g^{\prime}", hook] & G\text{-Hilb}_{\mathbb{C}}(\mathbb{P}^{4})\cong G\text{-Hilb}_{\mathbb{P}^{4}/G}(\mathbb{P}^{4}) \arrow[d, ""] \\
\mathbb{A}^{4}/G \arrow[r, "g", hook] & \mathbb{P}^{4}/G                                            
\end{tikzcd}$$
where $g$ and $g^{\prime}$ are open immersions. By Lemma \ref{Lemma for G-Hilbert of A4} and Proposition \ref{Proposition by Blume}, we have 
	$$G\text{-Hilb}_{\mathbb{C}}(\mathbb{A}^{4})\cong \mathbb{A}^{2}\times (\text{Bl}_{0}(\mathbb{A}^{2})/G)\cong \text{Bl}_{\mathbb{A}^{2}\times 0}(\mathbb{A}^{4})/G.$$
Taking $x_{i}=1$ for $i=0,1,2,3,4$ respectively, we get a Zariski cover of $G\text{-Hilb}_{\mathbb{C}}(\mathbb{P}^{4})$. Since the gluing of this cover is completely determined by the open embeddings
$\mathbb{A}^{4}\hookrightarrow \mathbb{P}^{4}$, one verifies readily that $G\text{-Hilb}_{\mathbb{C}}(\mathbb{P}^{4})\cong (\text{Bl}_{\mathbb{P}^{2}\coprod\mathbb{P}^{1}}(\mathbb{P}^{4}))/G$.
\end{proof}

	Let us give the following description of the $G-$Hilbert scheme $G\text{-Hilb}_{\mathbb{C}}(Y)$ for a non-Eckardt cubic threefold $(Y,\tau)$.
\begin{proposition}\label{G-Hilbert scheme of Y}
The $G$-Hilbert scheme $G\text{-Hilb}_{\mathbb{C}}(Y)\cong \text{Bl}_{\mathbb{P}^{2}\coprod\mathbb{P}^{1}}(\mathbb{P}^{4})/G\times_{\mathbb{P}^{4}/G} Y/G$  
\end{proposition}
\begin{proof}
Proceeding as in Lemma \ref{G-Hilbert scheme of P4}, we obtain the following Cartesian diagram:
$$\begin{tikzcd}
G\text{-Hilb}_{\mathbb{C}}(Y)\cong G\text{-Hilb}_{Y/G}(Y) \arrow[d, ""] \arrow[r, "f^{\prime}", hook] & G\text{-Hilb}_{\mathbb{C}}(\mathbb{P}^{4})\cong G\text{-Hilb}_{\mathbb{P}^{4}/G}(\mathbb{P}^{4}) \arrow[d, ""] \\
Y/G \arrow[r, "f", hook] & \mathbb{P}^{4}/G                                            
\end{tikzcd}$$
where $f$ and $f'$ are closed immersions. By Lemma \ref{G-Hilbert scheme of P4}, one has $G\text{-Hilb}_{\mathbb{C}}(\mathbb{P}^{4})\cong \text{Bl}_{\mathbb{P}^{2}\coprod\mathbb{P}^{1}}(\mathbb{P}^{4})$.
Hence, we obtain the isomorphism as desired
	$$G\text{-Hilb}_{\mathbb{C}}(Y)\cong \text{Bl}_{\mathbb{P}^{2}\coprod\mathbb{P}^{1}}(\mathbb{P}^{4})/G\times_{\mathbb{P}^{4}/G} Y/G.$$
\end{proof}
 
 	By \cite[Theorem 1.1]{BKR2001mckay}, the equivariant derived category $D^{b}_{G}(Y)$ is equivalent to the derived category $D^{b}(Z)$, where $Z$ is the irreducible component of $G\text{-Hilb}_{\mathbb{C}}(Y)$ containing the open subset of all reduced $G$-clusters. Note also that $Z$ is birational to $Y/G$. Under the isomorphism $G\text{-Hilb}_{\mathbb{C}}(Y)\cong \text{Bl}_{\mathbb{P}^{2}\coprod\mathbb{P}^{1}}(\mathbb{P}^{4})/G\times_{\mathbb{P}^{4}/G} Y/G$, the irreducible component $Z$ can be described as follows. 
\begin{proposition}\label{Theorem about the important component of G hilbert scheme}
	Suppose that the equation of $(Y,\tau)$ are given by Equations \eqref{eqn_nEkcardt_Y} and \eqref{eqn_nEcardt_involution}. Then the irreducible component $Z$ is isomorphic to $\text{Bl}_{C\coprod L}(Y)/G$ which is further isomorphic to $\text{Bl}_{C_{1}}(\mathbb{P}^{2}\times\mathbb{P}^{1})$ with $C_{1}$ being the curve in $\mathbb{P}^{2}\times\mathbb{P}^{1}$ cut out by $$V(g(x_{0},x_{1},x_{2}),x_{0}q_{0}(x_{3},x_{4})+x_{1}q_{1}(x_{3},x_{4})+x_{2}q_{2}(x_{3},x_{4}))\subset \mathbb{P}^{2}_{[x_{0}:x_{1}:x_{2}]}\times\mathbb{P}^{1}_{[x_{3}:x_{4}]}.$$ 
\end{proposition}
\begin{proof}
	The proof is quite similar with that of \cite[Proposition 4.1, Theorem 4.2]{Hu2023Equivariant}. By Proposition \ref{G-Hilbert scheme of Y}, the $G$-Hilbert scheme $G\text{-Hilb}_{\mathbb{C}}(Y)$ is isomorphic to $\text{Bl}_{\mathbb{P}^{2}\coprod\mathbb{P}^{1}}(\mathbb{P}^{4})/G\times_{\mathbb{P}^{4}/G} Y/G$. Since the only irreducible component of $\text{Bl}_{\mathbb{P}^{2}\coprod\mathbb{P}^{1}}(\mathbb{P}^{4})/G\times_{\mathbb{P}^{4}/G} Y/G$ which is isomorphic to $Y/G$ is $\text{Bl}_{C\coprod L}(Y)/G$, we have $Z=\text{Bl}_{C\coprod L}(Y)/G$. It then suffices to show that $\text{Bl}_{C\coprod L}(Y)/G\cong \text{Bl}_{C_{1}}(\mathbb{P}^{2}\times\mathbb{P}^{1})$. We shall verify this on an open affine cover. On one hand, the blowup $\text{Bl}_{\mathbb{P}^{2}\coprod \mathbb{P}^{1}}(\mathbb{P}^{4})$ is a subvariety of $\mathbb{P}^{4}\times \mathbb{P}^{2}\times \mathbb{P}^{1}$ with following equations:
\[
\begin{split}
& x_{0}'y_{1}'-x_{1}'y_{0}'=0;\
x_{1}'y_{2}'-x_{2}'y_{1}'=0;\
x_{0}'y_{2}'-x_{2}'y_{0}'=0;\
x_{3}'y_{4}'-x_{4}'y_{3}'=0
%& x_{0}'q_{0}+(x_{3}',x_{4}')+x_{1}'q_{1}(x_{3}',x_{4}')+x_{2}'q_{2}(x_{3}',x_{4}')+g(x_{0}',x_{1}',x_{2}')=0    
\end{split}
\]
where $[x_{0}':x_{1}':x_{2}':x_{3}':x_{4}']\times [y_{0}':y_{1}':y_{2}']\times [y_{3}':y_{4}']$ are the coordinates of $\mathbb{P}^{4}\times \mathbb{P}^{2}\times \mathbb{P}^{1}$. Taking the strict transform of $Y$ and then taking the quotient by $G$, we obtain a local description of the variety $\text{Bl}_{C\coprod L}(Y)/G$. On the other hand, the blowup $\text{Bl}_{C_{1}}(\mathbb{P}^{2}\times\mathbb{P}^{1})$ is cut out by the equation 
	$$y_0(x_{0}q_{0}(x_{3},x_{4})+x_{1}q_{1}(x_{3},x_{4})+x_{2}q_{2}(x_{3},x_{4}))+y_1g(x_{0},x_{1}, x_{2})=0$$
in $\mathbb{P}^{2}\times \mathbb{P}^{1}\times \mathbb{P}^{1}$ with coordinates $[x_{0}:x_{1}:x_{2}]\times [x_{3}:x_{4}]\times [y_{0}:y_{1}]$. Now letting $x_{0}'=y_{0}'=y_{3}'=1$, we obtain the following local equation of $\text{Bl}_{C\coprod L}(Y)/G$ in $\mathbb{A}^{4}_{(y_{1}',y_{2}',y_{4}',x_{3}'')}$ with $x_{3}'':=x_{3}'^{2}$
	$$x_{3}''(q_{0}(1,y_{4}')+x_{1}'q_{1}(1,y_{4}')+x_{2}'q_{2}(1,y_{4}'))+g(1,x_{1}',x_{2}')=0.$$
Taking $x_{0}=x_{3}=y_{1}=1$, the affine variety $\text{Bl}_{C_{1}}(\mathbb{P}^{2}\times \mathbb{P}^{1})|_{x_{0}=x_{3}=y_{1}=1}$ is given by $$V(y_{0}(q_{0}(1,x_{4})+x_{1}q_{1}(1,x_{4})+x_{2}q_{2}(1,x_{4}))+g(1,x_{1},x_{2}))\subset \mathbb{A}^{4}_{(x_{1},x_{2},x_{4},y_{0})}.$$ By mapping 
\[
\begin{split}
& x_{3}''\longrightarrow y_{0};\,\,
y_{1}'\longrightarrow x_{1}; \\
& y_{2}'\longrightarrow x_{2}; \,\,
  y_{4}'\longrightarrow x_{4},
\end{split}
\]
we get an isomorphism between the corresponding affine subvarieties. In a similar manner, we construct isomorphisms between other pairs of open affine varieties and verify that these local isomorphisms can be glued together. The proposition then follows.
\end{proof}

	Observe that the curve $C_1\subset\mathbb{P}^2\times \mathbb{P}^1$ admits a double covering map $p_1:C_1\to C$ to the smooth cubic curve $C$ in the fixed locus $Y^\tau=C\coprod L$ (project $\mathbb{P}^2\times \mathbb{P}^1$ to $\mathbb{P}^2$ and compare the equations of $C_1$ in Proposition \ref{Theorem about the important component of G hilbert scheme} and the equation of $C$ in Section \ref{section_cubic_involution} or \cite[\S2.2]{casalaina2022moduli}). The following proposition shows that $C_1$ is isomorphic to $\widetilde{C}$ which is also a double cover of the cubic curve $C\subset Y^\tau$.
	
\begin{proposition}\label{two_curves_constructed_same}
	The double cover $p_1:C_1\to C$ and the restriction of the discriminant double cover $\pi:\widetilde{C}\to C$ obtained by projecting $Y$ from the pointwise fixed line $L$ are isomorphic.
\end{proposition}
\begin{proof}
	Using the equations of $C_1$ in Proposition \ref{Theorem about the important component of G hilbert scheme}, it is not hard to show that the double covering map $p_1:C_1\to C$ is branched at the six intersection points $C\cap Q$ of the cubic component $C$ and the quadratic component $Q$ of the discriminant curve for the projection of $Y$ from the pointwise fixed line $L$. As a result, the curve $C_1$ is of (arithmetic) genus $4$. From singular Riemann-Hurwitz \cite[Equation 1.2]{garcia1996rational}, we deduce that the curve $C_{1}$ is smooth. As discussed in Section \ref{section_cubic_involution} (see also \cite[\S2.2]{casalaina2022moduli}), the double cover $\widetilde{C}\to C$ to the cubic curve $C$ is also branched at the six intersection points $C\cap Q$. Moreover, again by Proposition \ref{Theorem about the important component of G hilbert scheme} the fibers of $p_1:C_1\to C$ correspond to the $\tau$-invariant lines of $Y$ which are not pointwise fixed (more precisely, let $x\in C$ then the plane $\langle L,x\rangle\cap Y=L\cup l\cup l'$, and $p_1^{-1}(x)$ corresponds to the $\tau$-invariant lines $l$ and $l'$). By the construction of the restricted discriminant double cover $\pi:\widetilde{C}\to C$, we obtain a bijective morphism from $\widetilde{C}$ to $C_1$ which is an isomorphism.
\end{proof}	

    We conclude this section by completing the proof of Theorem \ref{thm_description_equivariant_derived_category}.

\begin{proof}[Proof of Theorem \ref{thm_description_equivariant_derived_category}]
 By Proposition \ref{Theorem about the important component of G hilbert scheme} and \cite[Theorem 1.1]{BKR2001mckay}, we have $D^b_{\mathbb{Z}_2}(Y)\cong D^b(Z)$ where $Z\cong Bl_{C_1}(\mathbb{P}^2\times\mathbb{P}^1)$. The result then follows from Proposition \ref{two_curves_constructed_same} and Orlov's blow-up formula \cite[Theorem 3.4]{beauville2016derived}.
\end{proof}

\section{Intermediate Jacobians of the equivariant Kuznetsov components}\label{section_equivariant_CT}
	In this section, we generalize the construction of the abstract intermediate Jacobian  in \cite{perry2020integral} for a smooth proper dg category. We then adapt this generalized construction to the equivariant Kuznetsov component $\Ku_{\mathbb{Z}_2}(Y)$ of a non-Eckardt cubic threefold $(Y,\tau)$, and give the first proof of Theorem \ref{main_theorem_first_equivariant_categorical_Torelli}. 

\subsection{Intermediate Jacobians of admissible subcategories}
	Let $X$ be a smooth projective variety. In this subsection, we very briefly review Perry's construction of the intermediate Jacobian of an admissible subcategory $\mathcal{A}$ of $\D^b(X)$. Let $\HH_{\bullet}(\mathcal{A})$ denote the Hochschild homology of $\mathcal{A}$, and let $\mathrm{HP}_{\bullet}(\mathcal{A})$ denote its periodic cyclic  homology. Also denoted by $K^{\mathrm{top}}_{\bullet}(\mathcal{A})$ the topological $\mathrm{K}$-group of $\mathcal{A}$ defined by Blanc in \cite[\S4.1]{blanc_2016}. Note that periodic cyclic homology and topological $\mathrm{K}$-theory are both $2$-periodic.

\begin{definition}[{\cite[Definition 5.24]{perry2020integral}}]\label{jacobian}
	Let $\mathcal{A}$ be an admissible subcategory of $\D^{b}(X)$. Consider the commutative diagram
 \[
\begin{tikzcd}
K_1^\mathrm{top}(\cA) \arrow[r, "\ch_1^{\mathrm{top}}"] \arrow[rrd, "\pi'"'] & \HP_1(\cA) \arrow[r, "\cong"] & \bigoplus_n \HH_{2n-1}(\cA) \arrow[d, "\pi"]   \\
                                                                           &                               & \HH_1(\cA) \oplus \HH_3(\cA) \oplus \cdots
\end{tikzcd}
 \]
with $\pi$ the natural projection and $\pi'$ the composition. The \emph{intermediate Jacobian of $\cA$} is the complex torus defined as follows
	$$J(\mathcal{A})= (\HH_{1}(\mathcal{A})\oplus \HH_{3}(\mathcal{A})\oplus\cdots)/\Gamma$$
where $\Gamma$ denotes the image of $\pi'$ and is a lattice.
\end{definition}

\begin{lemma}{(\cite[Lemma 3.9]{jacovskis2023infinitesimal})}\label{assumption}
 Let $X$ be a Fano threefold. Assume that the derived category $D^{b}(X)$ admits a semi-orthogonal decomposition $D^{b}(X)=\langle\Ku(X), E_{1}, E_{2}, \cdots, E_{n} \rangle$ where 
 $E_{1}, E_{2}, \ldots, E_{n}$ is an exceptional collection. Then we have $J(\Ku(X))\cong J(X)$.
\end{lemma}

\subsection{Intermediate Jacobians of smooth proper dg categories}\label{section_IJ_smooth_proper}
	Let $\cA$ be a smooth proper dg category. Since the non-commutative Hodge to de Rham spectral sequence of $\cA$ degenerates (see for instance \cite[Theorem 5.5]{kaledin2008non} and \cite[Theorem 5.4]{kaledin2017spectral}), we have the following Hodge filtration (where $\mathrm{HN}_{\bullet}(\mathcal{A})$ and $\HP_{\bullet}(\mathcal{A})$ denote the negative cyclic homology and the periodic cyclic homology of $\mathcal{A}$ respectively)
	$$\cdots\subset \mathrm{HN}_{-3}(\cA)\subset \mathrm{HN}_{-1}(\cA)\subset \mathrm{HN}_{1}(\cA)\subset \cdots \subset \mathrm{HP}_{1}(\cA).$$
In particular, there exists a short exact sequence
	$$0\rightarrow \mathrm{HN}_{2i-1}(\cA)\rightarrow  \mathrm{HN}_{2i+1}(\cA)\rightarrow \mathrm{HH}_{2i+1}(\cA)\rightarrow 0$$ for each $i\in \mathbb{Z}$. 

\begin{definition}\label{generalization_of_Alex_J}
   Let $\cA$ be a smooth proper dg category. We define the \emph{intermediate Jacobian} $J(\cA)$ of $\cA$ to be the following group
   $$J(\cA)=\frac{\mathrm{HP}_{1}(\cA)}{\mathrm{HN}_{-1}(\cA)+\mathrm{im}(K^{top}_{1}(\cA))}.$$
\end{definition}

\begin{remark}\label{remark_generalization_of_J_admissible}
	If $\cA$ is an admissible subcategory, then the group $J(\cA)$ is a complex torus which is isomorphic to the intermediate Jacobian of $\cA$ introduced in Definition \ref{jacobian}.
\end{remark}

	We also have the following lemma whose proof is straightforward. 
\begin{lemma}
	Let $\cA=\langle \cA_{1},\cA_{2}\rangle$ be a semi-orthogonal decomposition of a pre-triangulated dg category $\cA$. Then it holds that 
$$J(\cA)\cong J(\cA_{1})\oplus J(\cA_{2}).$$
Furthermore, if $\cA$, $\cA_{1}$ and $\cA_{2}$ are all admissible subcategories, then the above isomorphism is an isomorphism of complex tori.
\end{lemma}

\subsection{Intermediate Jacobians of the equivariant Kuznetsov components}
	Let $(Y,\tau)$ be a cubic threefold with a non-Eckardt type involution. Applying the intermediate Jacobian construction introduced in Definition \ref{generalization_of_Alex_J} to the equivariant Kuznetsov component $\Ku_{\mathbb{Z}_2}(Y)$	 (see Lemma \ref{cubic_threefold_equivariant_SOD}) which is a smooth dg category, we obtain the following proposition. 

\begin{proposition}\label{Lemma_intermediate_Jacobian_Ku_same_curve}
	Let $\widetilde{C}$ be the double cover of the cubic curve $C\subset Y^\tau$ as described in Section \ref{section_cubic_involution} (see also \cite[Proposition 2.5]{casalaina2022moduli}). Then we have an isomorphism of principally polarized abelian varieties
	$$J(\Ku_{\mathbb{Z}_2}(Y))\cong J(\widetilde{C}).$$    
\end{proposition}

\begin{proof}
	By Lemma \ref{cubic_threefold_equivariant_SOD}, the invariant derived category $D^b_{\mathbb{Z}_2}(Y)$ admits the following semi-orthogonal decomposition
	$$D^b_{\mathbb{Z}_2}(Y)=\langle\Ku_{\mathbb{Z}_2}(Y),\oh_Y,\oh_Y\otimes\chi_1,\oh_Y(1),\oh_Y(1)\otimes\chi_1\rangle.$$ 
By Theorem \ref{thm_description_equivariant_derived_category}, there exists another semi-orthogonal decomposition
	$$D^b_{\mathbb{Z}_2}(Y)=D^b(Z)=\langle D^b(\widetilde{C}),E_{ij}\rangle_{1\leq i\leq 3,1\leq j\leq 2}$$ 
where $E_{ij}$ with $1\leq i\leq 3$, $1\leq j\leq 2$ is an exceptional collection of line bundles on $\mathbb{P}^2\times\mathbb{P}^1$. 
We thus get 
	$$K^{top}_{1}(\Ku_{\mathbb{Z}_{2}}(Y))\cong K^{top}_{1}(D^{b}_{\mathbb{Z}_{2}}(Y))\cong K^{top}_{1}(D^{b}(\widetilde{C})\cong H^{1}(\widetilde{C},\mathbb{Z})$$
which is compatible with the Euler parings (cf.~\cite[Lemma 5.2]{perry2020integral}). Note that since the Euler paring on $K^{top}_{1}(D^{b}(\widetilde{C}))\cong H^{1}(C)$ coincides with the cohomology paring, the Euler paring on $K^{top}_{1}(\Ku_{\mathbb{Z}_{2}}(Y))$ is anti-symmetric. Now consider the following isomorphisms of the Hodge filtration on the dg categories. 
$$\xymatrix{\mathrm{HN}_{-1}(\Ku_{\mathbb{Z}_{2}}(Y))\ar[r]^{\hookrightarrow}\ar[d]_{\cong}&\mathrm{HP}_{1}(\Ku_{\mathbb{Z}_2}(Y)\ar[d]_{\cong}\\
\mathrm{HN}_{-1}(D^{b}(\widetilde{C}))\ar[r]^{\hookrightarrow}\ar[d]_{\cong}&\mathrm{HP}_{1}(D^{b}(\widetilde{C}))\ar[d]_{\cong}\\
H^{0,1}(\widetilde{C})\ar[r]^{\hookrightarrow}&H^{1}(\widetilde{C})}$$
The weight $1$ Hodge structure on $H^{1}(\widetilde{C},\mathbb{Z})$ induces a Hodge structure on $K^{top}_{1}(\Ku_{\mathbb{Z}_{2}}(Y))$ of weight $1$. Namely, 
$$K^{top}_{1}(\Ku_{\mathbb{Z}_{2}}(Y))\otimes \mathbb{C}\cong \mathrm{HP}_{1}(\Ku_{\mathbb{Z}_{2}}(Y))\cong \mathrm{HN}_{-1}(\Ku_{\mathbb{Z}_{2}}(Y))\oplus \overline{\mathrm{HN}_{-1}(\Ku_{\mathbb{Z}_{2}}(Y))}.$$
This shows that $J(\Ku_{\mathbb{Z}_{2}}(Y))$ constructed in Definition \ref{generalization_of_Alex_J} is a complex torus. Furthermore, the alternating Euler paring  on $K^{top}_{1}(\Ku_{\mathbb{Z}_{2}}(Y))$ induces a principal polarization on $J(\Ku_{\mathbb{Z}_{2}}(Y))$ making it a principally polarized abelian variety. Thus, we obtain the following isomorphism of principally polarized abelian varieties 
	$$J(\Ku_{\mathbb{Z}_2}(Y))\cong J(D^b_{\mathbb{Z}_2}(Y))\cong J(D^b(\widetilde{C}))\cong J(\widetilde{C}).$$
\end{proof}

	We are now ready to give the first proof or Theorem \ref{main_theorem_first_equivariant_categorical_Torelli}. 

\begin{theorem}[=Theorem \ref{main_theorem_first_equivariant_categorical_Torelli}]\label{theorem_equivariant_categorical_Torelli}
	Let $(Y,\tau)$ and $(Y',\tau')$ be general cubic threefolds with a non-Eckardt type involution. Assume that there is a Fourier-Mukai type equivalence $\Phi:\Ku_{\mathbb{Z}_2}(Y)\simeq\Ku_{\mathbb{Z}_2}(Y')$ between the equivariant Kuznetsov components. Then $(Y,\tau)\cong (Y',\tau')$.
\end{theorem}
\begin{proof}
By the assumption, there is an Fourier-Mukai equivalence $\Phi:\Ku_{\mathbb{Z}_2}(Y)\simeq\Ku_{\mathbb{Z}_2}(Y')$. Applying the construction in Definition \ref{generalization_of_Alex_J}, we obtain an isomorphism of principally polarized abelian varieties
	$$J(\Ku_{\mathbb{Z}_2}(Y))\cong J(\Ku_{\mathbb{Z}_2}(Y')).$$
Let $\widetilde{C}$ (respectively, $\widetilde{C'}$) be the double cover of the cubic component $C\subset Y^\tau$ (respectively $C'\subset (Y')^{\tau'}$) of the fixed locus $Y^\tau$ (respectively, $(Y')^{\tau'}$). By Proposition \ref{Lemma_intermediate_Jacobian_Ku_same_curve}, we get the following isomorphism of principally polarized abelian varieties
	$$J(\widetilde{C})\cong J(\widetilde{C'})$$
which further implies that $\widetilde{C}\cong \widetilde{C'}$. Since $(Y,\tau)$ and $(Y',\tau')$ are general non-Eckardt cubic threefolds, the restricted discriminant double covers $\widetilde{C}\to C$ and $\widetilde{C'}\to C'$, as well as the genus $4$ bielliptic curves $\widetilde{C}$ and $\widetilde{C'}$, are also general (cf.~\cite[Theorem 3.2]{casalaina2022moduli}). By the results in \cite{CDC05}, a general genus $4$ bielliptic curve admits a unique bielliptic structure (which can also be proved using a degeneration argument; for instance one considers an admissible double covering map from the union of three transverse elliptic curves $E_0\cup E_1\cup E_2$ with different, general $j$-invariants to the union of  three transverse lines $L_0\cup L_1\cup L_2$ in $\mathbb{P}^2$). As a result, one deduces that $\widetilde{C}\to C$ is isomorphic to $\widetilde{C'}\to C'$; this, together with \cite[Theorem 2.5]{casalaina2022moduli}, implies that the invariant parts of the intermediate Jacobians $J(Y)^{\tau}$ and $J(Y')^{\tau'}$ are isomorphic. Using \cite[Theorem 3.1]{casalaina2022moduli} we conclude that $(Y,\tau)\cong (Y',\tau')$. 

%Note that the Jacobian $J(\widetilde{C})$ is isogeneous to the product $J(C)\times P(\widetilde{C},C)$. Similarly, $J(\widetilde{C'})$ is isogeneous to $J(C')\times P(\widetilde{C'},C')$. Taking dual and applying \cite[Theorem 2.5]{casalaina2022moduli} (which shows that the invariant part $J(Y)^\tau$ is isomorphic to the dual abelian variety of the Prym variety $P(\widetilde{C},C)$ and the same holds for $J(Y')^{\tau'}$), one gets that $J(C)\times J(Y)^{\tau}$ is isogenous to $J(C')\times J(Y')^{\tau'}$. By \cite[Theorem 3.1, Theorem 3.4]{casalaina2022moduli}, $J(Y)^{\tau}$ and $J(Y')^{\tau'}$ are very general members in the moduli space $\mathcal{A}_3^{(1,2,2)}$. It follows that both of them are simple abelian varieties. From Poincar\'e reducibility theorem we deduce that $J(Y)^{\tau}$ is isogenous to $J(Y')^{\tau'}$ which further implies that $J(Y)^{\tau}\cong J(Y')^{\tau'}$ (taking abelian varieties up to isogeny of a given degree corresponds to taking the quotient of the Siegel upper half space by a larger arithmetic group of finite index;  taking a limit over all such quotients, one gets the assertion). Using \cite[Theorem 3.1]{casalaina2022moduli} we conclude that $(Y,\tau)\cong (Y',\tau')$.
\end{proof}

\section{Bridgeland moduli spaces on the equivariant Kuznetsov components }\label{section_reconstruct_Bridgeland_moduli_space}
	Let $(Y,\tau)$ be a cubic threefold with a non-Eckardt type involution. In this section, we reconstruct the fixed locus $F(Y)^{\tau}$ of the Fano surface of lines as a Bridgeland moduli space of stable objects in $\Ku_{\mathbb{Z}_2}(Y)$. As a consequence, we give the second proof of Theorem \ref{main_theorem_first_equivariant_categorical_Torelli}. 
	
	We refer the reader to \cite[Section 4.3]{perry2023moduli} for a review of the general theory of stability conditions on equivariant categories. In our situation, let $\sigma$ be a Serre-invariant stability condition on the Kuznetsov component $\Ku(Y)$. Since $\sigma$ is fixed by the non-Eckardt type involution $\tau$, it induces a stability condition $\sigma^{\mathbb{Z}_2}$ on the equivariant Kuznetsov component $\Ku_{\mathbb{Z}_2}(Y)$. Indeed, by \cite[Corollary 5.5]{pertusi2020some} and \cite[Theorem 4.20]{jacovskis2021categorical} there exists an element $\widetilde{g}\in\widetilde{GL}^+(2,\mathbb{R})$ such that $\sigma=\sigma_0\cdot \widetilde{g}$ where $\sigma_0$ is a tilted stability condition which is clearly fixed by $\tau$. It then follows that $\tau(\sigma)=\tau(\sigma_0\cdot \widetilde{g})=\sigma_0\cdot \widetilde{g}=\sigma$, and hence $\sigma$ is a $\mathbb{Z}_2$-fixed stability condition in the sense of \cite[Section 4.3]{perry2023moduli} and induces a stability condition $\sigma^{\mathbb{Z}_2}$ on $\Ku_{\mathbb{Z}_2}(Y)$. In other words, we have proved the following lemma.
	
\begin{lemma}\label{lem_existence_stability}
	Let $\Ku_{\mathbb{Z}_2}(Y)$ be the equivariant category of the Kuznetsov component $\Ku(Y)$ for a non-Eckardt cubic threefold $(Y,\tau)$ where $\mathbb{Z}_2=\langle 1,\tau\rangle$. Then there exists a stability condition $\sigma^{\mathbb{Z}_2}$ on $\Ku_{\mathbb{Z}_2}(Y)$. 
\end{lemma}

\begin{remark}\label{remk_existence_action}
	Besides Lemma \ref{cubic_threefold_equivariant_SOD}, one could also apply \cite[Corollary 4.10]{beckmann2023equivariant} to show that there exists a $\mathbb{Z}_2$-action on $\Ku(Y)$. Specifically, on one side, the Serre-invariant stability condition $\sigma$ is fixed by $\tau$, and by \cite[Corollary 6.11]{Li2024Higher} one gets $\Ku(Y)\simeq D^b(\mathcal{A})$ where $\mathcal{A}$ is the heart of $\sigma$. On the other side, there exists a $\tau$-invariant simple object $I_l$ (where $l\in F(Y)^{\tau}$ is a $\tau$-invariant line and $I_l$ denotes the ideal sheaf).
\end{remark}

	Let us recall the following theorem which will be useful.

\begin{theorem}{(\cite[Proposition 2.23]{polishchuk2006constant}, \cite[Theorem 1.1]{macri2007inducing}, \cite[Theorem 4.8]{perry2023moduli})} \label{thm_relation_stable_objects_equiva_original}
	In the above setup, let $\sigma:=(\cA,Z)$ be a $\mathbb{Z}_2$-fixed stability condition with respect to $v$. Define 
\begin{itemize}
    \item $v^{\mathbb{Z}_2}=v\circ\mathrm{Forg}_*: K_0(\Ku_{\mathbb{Z}_2}(Y))\rightarrow\mathcal{N}(\Ku(Y))$;
    \item $\cA^{\mathbb{Z}_2}=\{E\in\Ku_{\mathbb{Z}_2}(Y)\mid\mathrm{Forg}(E)\in\cA\}$;
    \item $Z^{\mathbb{Z}_2}=Z\circ v^{\mathbb{Z}_2}$
\end{itemize}
where $\mathrm{Forg}_*$ is the map on the Grothendieck groups induced by the forgetful functor. Then the pair $\sigma^{\mathbb{Z}_2}=(\cA^{\mathbb{Z}_2}, Z^{\mathbb{Z}_2})$ is a stability condition on $\Ku_{\mathbb{Z}_2}(Y)$ with respect to $v^{\mathbb{Z}_2}$. Moreover, for an object $E\in\Ku_{\mathbb{Z}_2}(Y)$ the following statements hold.
\begin{enumerate}
    \item The object $E$ is $\sigma^{\mathbb{Z}_2}$-semistable if and only if $\mathrm{Forg}(E)$ is $\sigma$-semistable. 
    \item If $E$ is $\sigma^{\mathbb{Z}_2}$-stable, then $\mathrm{Forg}(E)$ is $\sigma$-polystable. 
    \item If $\mathrm{Forg}(E)$ is $\sigma$-stable, then $E$ is $\sigma^{\mathbb{Z}_2}$-stable. 
\end{enumerate}
\end{theorem}

Let $\sigma$ be a Serre-invariant stability condition on $\Ku(Y)$. We now study the Bridgeland moduli space $\mathcal{M}_{\sigma^{\mathbb{Z}_2}}(\Ku_{\mathbb{Z}_2}(Y),\bv)$ of $\sigma^{\mathbb{Z}_2}$-stable objects in the equivariant Kuznetsov component $\Ku_{\mathbb{Z}_2}(Y)$. 

\begin{proposition}\label{invariant_lines_as_Bridgeland_moduli_space}
	Let $\sigma^{\mathbb{Z}_2}$ be the stability condition on $\Ku_{\mathbb{Z}_2}(Y)$ induced by a Serre-invariant stability condition $\sigma$ on $\Ku(Y)$. Then there is an isomorphism $$\mathcal{M}_{\sigma^{\mathbb{Z}_2}}(\Ku_{\mathbb{Z}_2}(Y),\bv)\cong F(Y)^{\tau}\coprod F(Y)^{\tau},$$ where $\bv=[I_l]$ denotes the class of the ideal sheaf of a ($\tau$-invariant) line $l\subset Y$.  
\end{proposition}
\begin{proof}
	Denote the trivial and the sign representations of $\mathbb{Z}_2$ by $\rho_0$ and $\rho_1$ respectively. First we show that $I_l\otimes \rho_i\in\Ku_{\mathbb{Z}_2}(Y)$ where $l$ is a $\tau$-invariant line in $Y$ (as discussed in Section \ref{section_cubic_involution}, either $l$ is the pointwise fixed line $L$ or it is parametrized by the double cover $\widetilde{C}$ of the cubic curve $C\subset Y^\tau$). Because $\tau^*(I_l)\cong I_l$, we have $I_l\in D^b_{\mathbb{Z}_2}(Y)$. By Lemma \ref{cubic_threefold_equivariant_SOD}, there exists a semi-orthogonal decomposition for $D^b_{\mathbb{Z}_2}(Y)$: 
	$$D^b_{\mathbb{Z}_2}(Y)=\langle\Ku_{\mathbb{Z}_2}(Y),\oh_Y,\oh_Y\otimes\chi_1,\oh_Y(1),\oh_Y(1)\otimes\chi_1\rangle.$$
Since $I_l\in\Ku(Y)$, it holds that $I_l\otimes \rho_0 \in\Ku_{\mathbb{Z}_2}(Y)$. Similarly, $I_l\otimes\rho_1\in\Ku_{\mathbb{Z}_2}(Y)$. From \cite[Lemma 5.13]{pertusi2020some} we deduce that $\mathrm{Forg}(I_l\otimes\rho_i)=I_l$ is $\sigma$-stable in $\Ku(Y)$. By Theorem \ref{thm_relation_stable_objects_equiva_original}, both $I_l\otimes\rho_0$ and $I_l\otimes\rho_1$ are $\sigma^{\mathbb{Z}_2}$-stable. Now let us show that any $\sigma^{\mathbb{Z}_2}$-stable object $E$ in $\Ku_{\mathbb{Z}_2}(Y)$ of class $\bv$ is isomorphic to either $I_l\otimes\rho_0$ or $I_l\otimes\rho_1$, where $l$ is a $\tau$-invariant line. It follows from Theorem \ref{thm_relation_stable_objects_equiva_original} that $\mathrm{Forg}(E)$ is $\sigma$-semistable in $\Ku(Y)$ of class $\bv$. Since $\bv$ is primitive, it is automatically $\sigma$-stable. By \cite[Lemma 5.13]{pertusi2020some}, we have $\mathrm{Forg}(E)\cong I_l$. Note that $\tau^*(E)\cong E$, and thus the underlying object of a $\sigma$-stable object $E\in\Ku_{\mathbb{Z}_2}(Y)$ is isomorphic to $I_l$ for $l\in F(Y)^{\tau}$. But there are two linearizations $E\rightarrow\tau^*(E)$: the trivial one $I_l\otimes\rho_0$ and the one $I_l\otimes\rho_1$ with a negative sign. Furthermore, they represent different points in the moduli space $\mathcal{M}_{\sigma}(\Ku_{\mathbb{Z}_2}(Y),\bv)$. We now conclude that $E$ is isomorphic to either $I_l\otimes\rho_0$ or $I_l\otimes\rho_1$, with $l\in F(Y)^\tau$. The proposition then follows.
\end{proof}

	In the remaining part of the section, we give the second proof of the main result Theorem \ref{main_theorem_first_equivariant_categorical_Torelli}  using Proposition \ref{invariant_lines_as_Bridgeland_moduli_space}. We will also need the following results. 

\begin{lemma}\label{equivalence_sendto_closed_point}
	Let $(Y,\tau)$ and $(Y',\tau')$ be cubic threefolds with a non-Eckardt type involution. Suppose that $\Phi:\Ku_{\mathbb{Z}_2}(Y)\simeq\Ku_{\mathbb{Z}_2}(Y')$ is an equivalence. Denote by $\mathcal{A}^{\mathbb{Z}_2}$ (respectively, $\mathcal{A}'^{\mathbb{Z}_2}$) the heart of the stability condition $\sigma^{\mathbb{Z}_2}$ (respectively, $(\sigma')^{\mathbb{Z}_2}$). Let $E\in\mathcal{A}^{\mathbb{Z}_2}$ be a $\sigma^{\mathbb{Z}_2}$-stable object of class $\bv=v^{\mathbb{Z}_2}(E)$. Then $\Phi(E)\in \mathcal{A}'^{\mathbb{Z}_2}$ and it is $(\sigma')^{\mathbb{Z}_2}$-stable. 
\end{lemma}
\begin{proof}
	WLOG, we may assume that $v^{\mathbb{Z}_2}(\Phi(E))=\bv$. Indeed, let $E$ be a $\sigma^{\mathbb{Z}_2}$-stable object of class $\bv=v^{\mathbb{Z}_2}(E)=1-l$. This implies that $v(\mathrm{Forg}(E))=\bv$, and hence by Proposition \ref{invariant_lines_as_Bridgeland_moduli_space} we get $\mathrm{Forg}(E)\cong E\cong I_l$ with $l\in F(Y)^{\tau}$. Now suppose that $v^{\mathbb{Z}_2}(\Phi(E))=v^{\mathbb{Z}_2}(\Phi(I_l))\neq \bv$. This means that 
$v(\mathrm{Forg}(\Phi(I_l)))\neq\bv$ but $v(\mathcal{S}^m_{\Ku(Y)}(\mathrm{Forg}(\Phi(I_l))))=\bv$ for some $m\in\mathbb{Z}$. Then it holds that $v^{\mathbb{Z}_2}(\mathcal{S}_{\Ku_{\mathbb{Z}_2}(Y)}^m\circ\Phi(I_l))=\bv$ where $\mathcal{S}_{\Ku_{\mathbb{Z}_2}(Y)}$ is the induced Serre functor on $\Ku_{\mathbb{Z}_2}(Y)$ (this functor exists since $\tau^*:\Ku(Y)\rightarrow\Ku(Y)$ is compatible with $\mathcal{S}_{\Ku(Y)}$). Next let us verify that $\Phi(E)\in\mathcal{A}'^{\mathbb{Z}_2}$. By the definition of the heart $\mathcal{A}'^{\mathbb{Z}_2}$ (see Theorem \ref{thm_relation_stable_objects_equiva_original}), it suffices to show that $\Phi(E)\in\Ku_{\mathbb{Z}_2}(Y')$ and $\mathrm{Forg}(\Phi(E))\in\mathcal{A}'$. Note that $\mathrm{Forg}(\Phi(E))\in\Ku(Y')$. Since $\mathrm{Forg}$ is a fully faithful functor and $\Phi$ is an equivalence, we have 
	$$\mathrm{Ext}^1(\mathrm{Forg}(\Phi(E)),\mathrm{Forg}(\Phi(E)))\cong\mathrm{Ext}^1(\Phi(E),\Phi(E))\cong\mathrm{Ext}^1(E,E)\cong\mathrm{Ext}^1(I_l,I_l)=\mathbb{C}^2.$$
As the homological dimension of the heart $\mathcal{A}'$ is $2$, the same argument as that of \cite[Lemma 6.6]{altavilla2019moduli} allows us to show that $\mathrm{Forg}(\Phi(E))\in\mathcal{A}'$. Finally we show that $\Phi(E)$ is $\sigma'^{\mathbb{Z}_2}$-stable. By Theorem \ref{thm_relation_stable_objects_equiva_original} we only need to verify that $\mathrm{Forg}(\Phi(E))\in\mathcal{A}'$ is $\sigma'$-stable. Again we have $\mathrm{Ext}^1(\mathrm{Forg}(\Phi(E)),\mathrm{Forg}(\Phi(E)))=\mathbb{C}^2$; Arguing as in the proof of \cite[Lemma 6.6]{altavilla2019moduli}, we get that $\mathrm{Forg}(\Phi(E))$ is $\sigma'$-stable which completes the proof. 
\end{proof}

	Applying Lemma \ref{equivalence_sendto_closed_point} to the equivalence $\Phi$ and its inverse, we obtain the following proposition.
\begin{proposition}\label{Equivalence_implies_bijection_closed_points}
	Notation as in Lemma \ref{equivalence_sendto_closed_point}. The equivalence $\Phi:\Ku_{\mathbb{Z}_2}(Y)\simeq\Ku_{\mathbb{Z}_2}(Y')$ induces a bijection map between the closed points of the moduli spaces 
$\mathcal{M}_{\sigma^{\mathbb{Z}_2}}(\Ku_{\mathbb{Z}_2}(Y),\bv)$ and $\mathcal{M}_{\sigma'^{\mathbb{Z}_2}}(\Ku_{\mathbb{Z}_2}(Y'),\bv)$.
\end{proposition}

	If we further assume that $\Phi:\Ku_{\mathbb{Z}_2}(Y)\simeq\Ku_{\mathbb{Z}_2}(Y')$ is a Fourier-Mukai type equivalence, then $\Phi$ induces an isomorphism between the Bridgeland moduli spaces.

\begin{proposition}\label{Fourier-Mukai_implies_isomorphism}
	Notation as in Lemma \ref{equivalence_sendto_closed_point}. Assume in addition that $\Phi:\Ku_{\mathbb{Z}_2}(Y)\simeq\Ku_{\mathbb{Z}_2}(Y')$ is a Fourier-Mukai type equivalence. Then $\Phi$ induces an isomorphism 
	$$\phi\colon\mathcal{M}_{\sigma^{\mathbb{Z}_2}}(\Ku_{\mathbb{Z}_2}(Y),\bv)\cong\mathcal{M}_{\sigma'^{\mathbb{Z}_2}}(\Ku_{\mathbb{Z}_2}(Y'),\bv).$$ 
\end{proposition}
\begin{proof}
As argued in the proof of Lemma \ref{equivalence_sendto_closed_point}, we assume that $v^{\mathbb{Z}_2}(\Phi(E))=\bv$. Since $\Phi$ is a Fourier-Mukai type equivalence and $\Phi$ induces the bijection of closed points between above two moduli spaces, the same argument as that of \cite[Theorem 8.3]{GLZ2021conics} allows us to obtain the required isomorphism. 
\end{proof}

	Now we are ready to give the second proof of Theorem \ref{main_theorem_first_equivariant_categorical_Torelli}.
\begin{theorem}[=Theorem \ref{main_theorem_first_equivariant_categorical_Torelli}]\label{alterntative_proof}
	Let $(Y,\tau)$ and $(Y',\tau')$ be general cubic threefolds with a non-Eckardt type involution. Assume that there is a Fourier-Mukai type equivalence $\Phi:\Ku_{\mathbb{Z}_2}(Y)\simeq\Ku_{\mathbb{Z}_2}(Y')$ between the equivariant Kuznetsov components. Then $(Y,\tau)\cong (Y',\tau')$.
\end{theorem}
\begin{proof}
By Proposition \ref{Fourier-Mukai_implies_isomorphism}, we have the isomorphism 
$$\phi\colon \mathcal{M}_{\sigma^{\mathbb{Z}_2}}(\Ku_{\mathbb{Z}_2}(Y),\bv)\cong\mathcal{M}_{\sigma'^{\mathbb{Z}_2}}(\Ku_{\mathbb{Z}_2}(Y'),\bv).$$
Using Proposition \ref{invariant_lines_as_Bridgeland_moduli_space}, we get the isomorphism
$$\phi:F(Y)^{\tau}\coprod F(Y)^{\tau}\rightarrow F(Y')^{\tau'}\coprod F(Y')^{\tau'}.$$
As recalled in Section \ref{section_cubic_involution} (see also \cite[\S1.3]{casalaina2022moduli}), the fixed locus $F(Y)^{\tau}$ of Fano surface of lines in $Y$ consists of two connected components: a point corresponding to the pointwise fixed line $L\subset Y$ and a smooth curve $\widetilde{C}$ of genus $4$ which is the double cover of the cubic curve $C\subset Y^\tau$ and parameterizes other $\tau$-invariant lines $l\subset Y$. Similar statements hold for $F(Y')^{\tau'}$.
It can then be readily shown that the isomorphism $\phi:F(Y)^{\tau}\coprod F(Y)^{\tau}\rightarrow F(Y')^{\tau'}\coprod F(Y')^{\tau'}.$ gives rise to an isomorphism between $\widetilde{C}$ and $\widetilde{C'}$. Arguing as in the proof of Theorem \ref{theorem_equivariant_categorical_Torelli}, we conclude that $(Y,\tau)\cong (Y',\tau')$.
\end{proof}

\section{Chow theory of the equivariant Kuznetsov components}\label{section_chow_theory_Kuznetsov_components}
	In this section let us discuss the third proof of Theorem \ref{main_theorem_first_equivariant_categorical_Torelli} for very general cubic threefolds admitting a non-Eckardt type involution. Specifically, with notation remaining the same as in Theorem \ref{main_theorem_first_equivariant_categorical_Torelli}, we will show that the Fourier-Mukai equivalence $\Phi$ of the equivariant Kuznetsov components induces an isomorphism between the invariant parts $A^{2}_{\mathbb{Z}_{2}, \mathbb{Q}}(Y)$ and $A^{2}_{\mathbb{Z}_{2}, \mathbb{Q}}(Y')$ of the groups of algebraically trivial cycles. The proof can then be completed using this result and the Abel--Jacobi maps.

	We begin by giving the general setup of Fourier-Mukai functors in the equivariant setting (which will be used in the remaining part of this section). Let $G$ and $H$ be finite groups. Suppose that $X$ is a smooth variety with a $G$-action and let $Y$ be a smooth variety with an action by $H$. Let us also equip $Y$ with the trivial $G$-action so that $Y$ admits an action by $G\times H$.
\begin{definition}(Fourier-Mukai functors)
	Let $P\in D_{G\times H}^{b}(X\times Y)$. Denote the projections by $p_1$ and $p_2$. 
	$$\xymatrix{&X\times Y\ar[dl]_{p_{1}}\ar[dr]^{p_{2}}&\\X&&Y}$$	
The \emph{Fourier-Mukai functor} associated with $P$ is defined as
	$$\Phi_{P}: D_{G}^{b}(X) \longrightarrow D_{H}^{b}(Y), \quad E\mapsto (\mathbb{R}p_{2}(\mathbb{L}p^{\ast}_{1}(E)\otimes^{\mathbb{L}} P))^{\mathbb{R}G}$$
 where $\mathbb{R}G$ denotes the derived functor of taking invariant sections $(\cdot)^{G}$.
\end{definition}

\begin{definition}
	Let $\mathcal{A}$ (respectively, $\mathcal{B}$) be an admissible subcategory of $D_{G}^{b}(X)$ (respectively, $D_{H}^{b}(Y)$). A functor $\Phi: \mathcal{A}\longrightarrow \mathcal{B}$ is \emph{of Fourier-Mukai type} if there exists $P\in D_{G\times H}^{b}(X\times Y)$ such that the following diagram is commutative.
 $$
 \xymatrix{D_{G}^{b}(X)\ar[d]^{pr}\ar[r]^{\Phi_{P}}&D_{H}^{b}(Y)\\
 \mathcal{A}\ar[r]^{\Phi}&\mathcal{B}\ar[u]
 }
 $$
\end{definition}

	For the definitions of the derived functors $Ind^{G\times G}_{G}$ and $Res^{G\times G}_{G}$ used below, we refer the reader to \cite[Section 2]{BFK}. 
\begin{lemma}
 The Fourier-Mukai functor with kernel $Ind^{G\times G}_{G}\Delta_{\ast}\mathcal{O}_{X}\in ~D^{b}_{G\times G}(X\times X)$ is the identity functor on $D^{b}_{G}(X)$.  
\end{lemma}
\begin{proof}
 Let $E\in D^{b}_{G}(X)$, then it holds that
 \begin{align*}
 (Rp_{\ast,2}(p^{\ast}_{1}E\otimes Ind^{G\times G}_{G}\Delta_{\ast}\mathcal{O}_{X}))^{\mathbb{R}G}
 	\cong \,& (p_{2,\ast}Ind^{G\times G}_{G}(Res^{G\otimes G}_{G}p^{\ast}_{1}E\otimes \Delta_{\ast}\mathcal{O}_{X}))^{\mathbb{R}G}\\
	\cong \,&(p_{\ast,2}Ind^{G\times G}_{G}\Delta_{\ast}(\Delta^{\ast}Res^{G\times G}_{G}p^{\ast}_{1}E\otimes \mathcal{O}_{X}))^{\mathbb{R}G}\\
  	\cong \,& (p_{\ast,2}Ind^{G\times G}_{G}\Delta_{\ast}E)^{\mathbb{R}G}\\
  	\cong \,& E.
\end{align*}
Specifically, the first isomorphism follows from the projection formula for $Ind^{G\times G}_{G}$ and $Res^{G\times G}_{G}$ (e.g., \cite[Lem.~2.16 (c)]{BFK}). The second isomorphism follows from the projection formula for $\Delta^{\ast}$ and $\Delta_{\ast}$. Since $(p_{\ast,2}Ind^{G\times G}_{G}\Delta_{\ast})^{\mathbb{R}G}$ is the right adjoint to $\Delta^{\ast}Res^{G\times G}_{G}p^{\ast}_{2}\cong Id$, we get 
	$$(p_{\ast,2}Ind^{G\times G}_{G}\Delta_{\ast})^{\mathbb{R}G}\cong Id.$$
The last isomorphism then follows.
\end{proof}

	Suppose that there is a Fourier--Mukai functor $\Phi :\mathcal A \to \mathcal B$ between the admissible subcategories of $D_{G}^{b}(X)$ and $D_{H}^{b}(Y)$. Then $\Phi$ induces maps between the invariant Grothendieck groups, the invariant Chow groups and the invariant parts of the groups of algebraically trivial cycles (with coefficients in $\mathbb{Q}$). In other words, we have the following commutative diagram. (Note that besides taking Mukai vectors $\nu(-)$ one could also Chern classes and Chern characters (as in \cite[\S8]{kuznetsovnonclodedfield2021}) which also give maps $K_G(X)\to \operatorname{CH}_{G,\mathbb Q}^\bullet(X)$ and $K_G(Y)\to \operatorname{CH}_{G,\mathbb Q}^\bullet(Y)$.)	
$$
\xymatrix{
 \mathcal A \ar[r] \ar@{->}@/^2pc/[rrr]^\Phi& D_{G}^{b}(X)\ar[d]^{[-]}\ar[r]^{\Phi_{P}}&D_{H}^{b}(Y) \ar[r] \ar[d]^{[-]}& \mathcal B\\
 & K_G(X)\ar[r]^{\Phi_{[P]}} \ar[d]^{\nu(-)}&K_H(Y) \ar[d]^{\nu(-)} & \\
 &\operatorname{CH}_{G,\mathbb Q}^\bullet(X) \ar[r]^{\Phi_{\nu(P)}} & \operatorname{CH}_{H,\mathbb Q}^\bullet(Y)&  \\
 &\operatorname{A}_{G,\mathbb Q}^\bullet(X) \ar[r]^{\Phi_{\nu(P)}} \ar@{^(->}[u]& \operatorname{A}_{H,\mathbb Q}^\bullet(Y) \ar@{^(->}[u] &  \\
 }
$$

	In the case when $X$ and $Y$ are rationally connected threefolds, and $\mathcal A$ and $\mathcal B$ are components of semi-orthogonal decompositions, the techniques in \cite[\S8]{kuznetsovnonclodedfield2021} allows us to prove the theorem.

\begin{theorem}\label{thm_general_result_chow_isom}
	Let $X$ and $Y$ be smooth projective rationally connected threefolds, both with the action of a finite group $G$. Suppose that there are semi-orthogonal decompositions 
\begin{align*}
D^b(X)&=\langle \mathcal A,{}^\perp\mathcal{A}\rangle\\
D^b(Y)&=\langle \mathcal{B},{}^{\perp}\mathcal{B}\rangle
\end{align*}
such that ${}^{\perp}\mathcal{A}$ and ${}^{\perp}\mathcal{B}$ are generated by exceptional collections. If there exists a Fourier-Mukai type equivalence \[\Phi:\mathcal A_G\stackrel{\cong}{\longrightarrow} \mathcal B_G\] between the equivariant components, then we have an isomorphism between the invariant parts of the groups of algebraically trivial cycles with rational coefficients
	$$\xymatrix{A^{2}_{G, \mathbb{Q}}(X)\ar@/_/[r]_{\frac{c_{3}(P)}{2}}& A^{2}_{G, \mathbb{Q}}(Y)\ar@/_/[l]_{\frac{c_{3}(P')}{2}}}.$$

    Consequently, there is an isogeny 
$$
J_G(X)\to J_G(Y)
$$
between the $G$-invariant parts of the intermediate Jacobians.
\end{theorem}
\begin{proof}
Let $P$ be the kernel of the Fourier-Mukai type equivalence $\Phi$. We first show that $\Phi_{P}$ defines an isomorphism of the groups of algebraically trivial cycles $\Phi_{v(P)}: A^{\ast}_{G, \mathbb{Q}}(X)\rightarrow A^{\ast}_{G, \mathbb{Q}}(Y)$ where $v(P):=\mathrm{ch}(P)\sqrt{\mathrm{td}(X\times Y)}$. Since $X$ and $Y$ are rationally connected threefolds, we have $\Pic^0=0$ and hence the only non-trivial part of $A^*_\mathbb{Q}$ is $A^2_\mathbb{Q}$. Write $j: \cA_G \hookrightarrow D^{b}_{G}(X)$ as the embedding functor whose left adjoint $pr$ is the projection functor. Let $\Phi^{-1}$ denote the inverse of $\Phi$. Then $\Phi_{P^{'}}=j \circ\Phi^{-1}\circ pr$ is the left adjoint of $\Phi_{P}$ where $P'=(P^{T})^{\vee}\otimes p^{\ast}_{2}\omega_{Y}[3]$. Similar to \cite[Proposition 5.1]{caldararu2003mukai} and \cite[\S3]{Caldararu2010}, we have a natural morphism $P'\circ P\rightarrow Ind^{G\times G}_{G}\Delta_{\ast}\mathcal{O}_{X}$ corresponding to the natural transformation $\Phi_{P'}\circ \Phi_{P}\implies Id_{D^{b}_{G}(X)}$. There then exists a triangle
  $$P'\circ P\rightarrow Ind^{G\times G}_{G}\Delta_{\ast}\mathcal{O}_{X}\rightarrow \mathcal{Q}$$
 where $\mathcal{Q}$ is generated by $E_{i}\boxtimes E_{j}$ with $E_{i}$ and $E_{j}$ exceptional objects (see also Section \ref{section_FM_functors}). It follows that 
 	$$\Phi_{v(P')}\circ\Phi_{v(P)}=\Phi_{v(P'\circ P)}=\Phi_{Ind^{G\times G}_{G}\Delta_{\ast}\mathcal{O}_{X}}-\Phi_{v(\mathcal{Q})}.$$
According to \cite[Lemma 8.2(iii), Lemma 8.3(iii)]{kuznetsovnonclodedfield2021}, we have $\Phi_{v(\bullet)}=\Phi_{ch_{3}(\bullet)}=\Phi_{\frac{c_{3}(\bullet)}{2}}$. From \cite[Lemma 8.2(iii)]{kuznetsovnonclodedfield2021} we deduce that $\Phi_{v(\mathcal{Q})}=0$. As a result, $\frac{c_{3}(P')}{2}\circ \frac{c_{3}(P)}{2}$ is an isomorphism between $\mathbb{A}^{2}_{G,\mathbb{Q}}(X)$ and $\mathbb{A}^{2}_{G,\mathbb{Q}}(Y)$. A similar argument for the natural transformation $\Phi_{P}\circ\Phi_{P'}\implies Id_{D^{b}_{G}(Y)}$ shows that $\frac{c_{3}(P)}{2}\circ \frac{c_{3}(P')}{2}$ is an isomorphism from $\mathbb{A}^{2}_{G,\mathbb{Q}}(Y)$ to $\mathbb{A}^{2}_{G,\mathbb{Q}}(X)$. Thus, we obtain the following isomorphisms
$$\xymatrix{A^{2}_{G, \mathbb{Q}}(X)\ar@/_/[r]_{\frac{c_{3}(P)}{2}}& A^{2}_{G, \mathbb{Q}}(Y)\ar@/_/[l]_{\frac{c_{3}(P')}{2}}}.$$
In fact, one can see in the argument above that we only needed to invert $2$, and so the entire discussion above holds with $\mathbb Z[1/2]$ coefficients, not just with $\mathbb Q$ coefficients.

Next, let us show that there is an isogeny $J_G(X)\rightarrow J_G(Y)$. Note that we have the following commutative diagram.
\begin{equation}\label{E:A2GJG}
\xymatrix{
A^2_G(X)\ar[r] \ar@{-->>}[d] & A^2(X) \ar@{->>}[d]\\
J_G(X)\ar[r] & J(X)
}    
\end{equation}
where the horizontal morphisms are injective. The composition $A^2_G(X)\hookrightarrow A^2(X)\twoheadrightarrow J(X)$ of the inclusion map and the Abel--Jacobi map lands inside the invariant part $J_G(X)$ of intermediate Jacobian, as the Abel--Jacobi map is equivariant with respect to the $G$-action. Note that the horizontal inclusions in \eqref{E:A2GJG} split up to inverting $|G|$, and this therefore gives the vertical surjection on the left, as both the source and target are divisible groups.  

%Moreover, the kernel $P$ of Fourier-Mukai type equivalence $\Phi:\cA_G\simeq\cB_G$ induces the morphism $A^2(X)\rightarrow A^2(Y)$ such that it preserves the invariant parts. By the universal property of the algebraic representative (recall that the intermediate Jacobian for codimension-2 cycles is the algebraic representative), we get the diagram below. 

Any $G$-invariant correspondence induces a morphism $A^2(X)\rightarrow A^2(Y)$ that preserves the invariant parts. We get the diagram below by the universal property of the algebraic representative (recall that the intermediate Jacobian for codimension-2 cycles is the algebraic representative). 
$$
\xymatrix{
A^2(X) \ar[r] \ar@{->>}[d]& A^2(Y) \ar@{->>}[d]\\
J(X) \ar[r] & J(Y)
}
$$
Now we can use \eqref{E:A2GJG}  and any equivariant correspondence 
%applied to $X$ and the splitting of \eqref{E:A2GJG} applied to $Y$, after inverting $|G|$, 
to obtain the following commutative diagram.
$$
\xymatrix@R=1em{
A^2_{G}(X)\ar@{^(->}[rd] \ar@{->>}[ddd] \ar@{->}[rrr]&&&A^2_{G}(Y)\ar@{->>}[ddd] \ar@{_(->}[ld]\\
 & A^2_{}(X) \ar@{->>}[d] \ar[r]& A^2_{}(Y)  \ar@{->>}[d]& \\
& J(X)_{}  \ar[r] & J(Y)_{}  &\\
J_G(X)_{}\ar@{^(->}[ru] \ar@{-->}[rrr]&&& J_G(Y)_{} \ar@{_(->}[lu]
}
$$
This is all functorial, and takes the invariant parts to the invariant parts, giving rise to the dashed arrow in the diagram above. %The middle vertical arrows are isomorphisms. 

    Tensoring the above diagram by $\mathbb Z[1/2]$, and then 
applying this to the correspondences $\frac{c_3(P)}{2}$ and $\frac{c_3(P)'}{2}$ inducing the isomorphism of the groups $A^2_{G, \mathbb{Z}[1/2]}(X)\stackrel{\cong}{\to} A^2_{G, \mathbb{Z}[1/2]}(Y)$, together with a diagram chase using the surjectivity of the vertical arrows on the outside of the diagram above, gives the isomorphism of the $G$-invariant parts of the intermediate Jacobians, in the category of $\mathbb Z[1/2]$-Hodge structures (i.e., up to isogenies of orders a power of $2$).
\end{proof}

\begin{remark}
    The proof of the theorem actually shows that the kernel of the isogeny of $J_G(X)$ and $J_G(Y)$ can have the order of divisible at most by a power of $2$.
\end{remark}

\begin{remark} 
    In Theorem \ref{thm_general_result_chow_isom}, if $G$ is the trivial group, then the theorem only implies that $J(X)$ is isogenous to $J(Y)$.  However,  in this case, one can show that $J(X)$ is isomorphic to $J(Y)$ using topological K-groups $K_1^{top}$. 
\end{remark}

    We have the following result to the case when $(Y,\tau)$ and $(Y',\tau')$ are cubic threefolds with a non-Eckardt type involution.
\begin{corollary}\label{chowgroupisomorphism}
	Let $(Y,\tau)$ and $(Y',\tau')$ be cubic threefolds with non-Eckardt type involutions. If there exists a Fourier-Mukai type equivalence \[\Phi:\Ku_{\mathbb{Z}_2}(Y)\stackrel{\cong}{\to} \Ku_{\mathbb{Z}_2}(Y')\] between the equivariant Kuznetsov components, then we have an isomorphism
	$$\xymatrix{A^{2}_{\mathbb{Z}_{2}, \mathbb{Q}}(Y)\ar@/_/[r]_{\frac{c_{3}(P)}{2}}& A^{2}_{\mathbb{Z}_{2}, \mathbb{Q}}(Y')\ar@/_/[l]_{\frac{c_{3}(P')}{2}}}.$$ 
Furthermore, the invariant parts $J(Y)^{\tau}$ and $J(Y')^{\tau}$ of the intermediate Jacobians of cubic threefolds $Y$ and $Y'$ are isogenous. 
\end{corollary}

\begin{remark}
	Note that up to this point in this section, all the results stated so far hold over any field, after one replaces intermediate Jacobians with algebraic representatives. Over non-closed fields, one only needs to require that $X$ and $Y$ be geometrically rationally connected. See \cite{ACMV19} for similar arguments over perfect fields, and \cite{ACMV23}, which explains how to remove the hypothesis that the field is perfect. 
\end{remark}

	Let us conclude this section by giving the third proof of a slightly weaker version of Theorem \ref{main_theorem_first_equivariant_categorical_Torelli} using Theorem \ref{thm_general_result_chow_isom} and Corollary \ref{chowgroupisomorphism}.
\begin{theorem}[Theorem \ref{main_theorem_first_equivariant_categorical_Torelli} for very general non-Eckardt cubic threefolds]\label{alterntative_proof_chow}
	Let $(Y,\tau)$ and $(Y',\tau')$ be very general cubic threefolds with a non-Eckardt type involution. Assume that there is a Fourier-Mukai type equivalence $\Phi:\Ku_{\mathbb{Z}_2}(Y)\simeq\Ku_{\mathbb{Z}_2}(Y')$ between the equivariant Kuznetsov components. Then $(Y,\tau)\cong (Y',\tau')$.
\end{theorem}
\begin{proof}
	By Corollary \ref{chowgroupisomorphism}, we have $A^{2}_{\mathbb{Z}_{2}, \mathbb{Q}}(Y)\cong A^{2}_{\mathbb{Z}_{2}, \mathbb{Q}}(Y')$. Composing with the Abel--Jacobi maps, we get $J(Y)^{\tau}\cong J(Y')^{\tau}$ up to isogeny. By \cite[Theorem 3.1, Theorem 3.4]{casalaina2022moduli}, $J(Y)^{\tau}$ and $J(Y')^{\tau'}$ are very general members in the moduli space $\mathcal{A}_3^{(1,2,2)}$. It then follows that $J(Y)^{\tau}\cong J(Y')^{\tau'}$ (taking abelian varieties up to isogeny of a given degree corresponds to taking the quotient of the Siegel upper half space by a larger arithmetic group of finite index; taking a limit over all such quotients, one gets the assertion). By \cite[Theorem 3.1]{casalaina2022moduli} we obtain that $(Y,\tau)\cong (Y',\tau')$.	
\end{proof}

\section{An equivariant infinitesimal categorical Torelli theorem}\label{section_inf_cat_Torelli}
	In this section, we first prove Theorem \ref{theorem_second_main_result} on the compatibility of the algebra structures of the Hochschild cohomology and of the module structures of the Hochschild homology over the Hochschild cohomology (see Section \ref{sec_compactibility_of_additivity}). We then show that the action of an automorphism of a smooth projective variety on Hochschild cohomology is compatible with the twisted HKR isomorphism (see Section \ref{sec_equivariant_HKR}) which generalizes a result in \cite{macri2009infinitesimal}. Putting these together, we give a proof of Theorem \ref{main_theorem_comm_diag_equivariant} (which can be viewed as an equivariant version of the infinitesimal categorical Torelli theorem) in Section \ref{subsection_equivariant_inf_categorical_Torelli}. Most results in this section hold in a more general setting (not just for cubic threefolds with a non-Eckardt type involution); we hope they are of independent interest. 

\subsection{Compatibility of algebra and module structures on Hochschild (co)homology}\label{sec_compactibility_of_additivity}

\subsubsection{Functors of Fourier-Mukai type} \label{section_FM_functors}
	In this subsection, we define Fourier-Mukai type functors between semi-orthogonal components and discuss some of the basic properties from the kernel perspective; these partially generalize some of the results in \cite{huybrechts2016hochschild}. Specifically, let $X$ be a smooth projective variety. Suppose that there is a semi-orthogonal decomposition  
	$$D^{b}(X)=\Big\langle \mathcal A_{1}, \mathcal{A}_{2}, \cdots, \A_{m}\Big\rangle.$$
Then $D^{b}(X\times X)$ admits a semi-orthogonal decomposition
	$$D^{b}(X\times X)=\Big\langle\mathcal{A}_{1X},\cdots,\mathcal{A}_{mX} \Big\rangle$$
with $\mathcal{A}_{iX}=\mathcal{A}_{i}\boxtimes D^{b}(X)$. Let us write $j'_{i}$ as the embedding $\mathcal{A}_{iX}\hookrightarrow D^{b}(X\times X)$ and write $j'^{\ast}_{i}$ as the associated projection functor. Note that $j'^{\ast}_{1}$ is the left adjoint of $j'_{1}$. Let 
 	$$\mathcal{B}_{i}=^{\perp}\Big\langle \mathcal{A}_{1},\cdots,\mathcal{A}_{i-1},\mathcal{A}_{i+1},\cdots, \mathcal{A}_{m} \Big\rangle.$$
It holds that $\mathcal{B}_{i}^{\vee}=\mathrm{R}\mathcal{H}om(\mathcal{B}_{i},\mathcal{O}_{X})$. There is also a semi-orthogonal decomposition
	$$D^{b}(X)=\Big\langle\mathcal{B}^{\vee}_{1},\mathcal{B}^{\vee}_{2},\cdots,\mathcal{B}^{\vee}_{m}\Big\rangle.$$
Thus we obtain the following semi-orthogonal decomposition of $D^{b}(X\times X)$
	$$D^{b}(X\times X)=\Big\langle \mathcal{A}_{1}\boxtimes \mathcal{B}^{\vee}_{1},\cdots, \mathcal{A}_{1}\boxtimes\mathcal{B}^{\vee}_{m}, \mathcal{A}_{2}\boxtimes \mathcal{B}^{\vee}_{1},\cdots, \mathcal{A}_{2}\boxtimes \mathcal{B}^{\vee}_{m},\cdots, \mathcal{A}_{m}\boxtimes \mathcal{B}^{\vee}_{1}, \cdots, \mathcal{A}_{m}\boxtimes \mathcal{B}^{\vee}_{m}\Big\rangle$$
where $\mathcal{A}_{i}\boxtimes \mathcal{B}^{\vee}_{j}$ denotes the minimal sub-triangulated category containing objects $E_{i}\boxtimes F_{j}$ with $E_{i}\in \mathcal{A}_{i}$ and $F_{j}\in \mathcal{B}^{\vee}_{j}$. Write $j_{i}: \mathcal{A}_{i}\boxtimes \mathcal{B}^{\vee}_{i}\hookrightarrow D^{b}(X\times X)$ and let $j^{\ast}_{i}$ be the projection functor to $\mathcal{A}_{i}\boxtimes \mathcal{B}^{\vee}_{i}$; note that $j^{\ast}_{1}$ is the left adjoint of $j_{1}$. For the later use, we also define $P_{i}:=j'^{\ast}_{i}\Delta_{\ast}\mathcal{O}_{X}$. According to \cite[Proposition 3.8]{kuznetsov2009hochschild}, one has $P_{i}\cong j^{\ast}\Delta_{\ast}\mathcal{O}_{X}$.
\begin{lemma}\label{intersection}
Notation as above. It holds that $\mathcal{A}_{i}\boxtimes \mathcal{B}_{j}^{\vee}=\mathcal{A}_{iX}\cap_{X}\mathcal{B}^{\vee}_{j}$ for $1\leq i,j\leq m$.  
\end{lemma}
\begin{proof}
	On one hand, it is clear that $\mathcal{A}_{i}\boxtimes \mathcal{B}_{j}^{\vee}\subset\mathcal{A}_{iX}\cap_{X}\mathcal{B}^{\vee}_{j}$. On the other hand, for any object $E\in \mathcal{A}_{s}\boxtimes \mathcal{B}_{t}^{\vee}$ and $F\in \mathcal{A}_{iX}\cap_{X}\mathcal{B}^{\vee}_{j}$, we have $\Hom(F,E)=0$ whenever $t<j$ or $s<i$:
	$$\Hom(\mathcal{A}_{iX}\cap_{X}\mathcal{B}^{\vee}_{j},\mathcal{A}_{s}\boxtimes \mathcal{B}^{\vee}_{t})=0,\quad t< j\ \text{or}\ s< i.$$
Similarly, we also have    
	$$\Hom(\mathcal{A}_{s}\boxtimes \mathcal{B}^{\vee}_{t},\mathcal{A}_{iX}\cap_{X}\mathcal{B}^{\vee}_{j})=0, \quad t> j\ \text{or}\ s> i.$$
Since there exists a semi-orthogonal decomposition
	$$D^{b}(X\times X)=\Big\langle \mathcal{A}_{1}\boxtimes \mathcal{B}^{\vee}_{1},\cdots, \mathcal{A}_{1}\boxtimes\mathcal{B}^{\vee}_{m}, \mathcal{A}_{2}\boxtimes \mathcal{B}^{\vee}_{1},\cdots, \mathcal{A}_{2}\boxtimes \mathcal{B}^{\vee}_{m},\cdots, \mathcal{A}_{m}\boxtimes \mathcal{B}^{\vee}_{1}, \cdots, \mathcal{A}_{m}\boxtimes \mathcal{B}^{\vee}_{m}\Big\rangle,$$
it holds that 
	$$\mathcal{A}_{iX}\cap_{X}\mathcal{B}^{\vee}_{j}\subset \mathcal{A}_{i}\boxtimes \mathcal{B}_{j}^{\vee}.$$
Thus $\mathcal{A}_{i}\boxtimes \mathcal{B}_{j}^{\vee}=\mathcal{A}_{iX}\cap_{X}\mathcal{B}^{\vee}_{j}$.
\end{proof}
  
	Now consider the following diagram.
	$$\xymatrix{&X\times X\ar[dr]^{p_{2}}\ar[dl]_{p_{1}}&\\ X&&X}$$
For any object $E\in D^{b}(X\times X)$, the \emph{Fourier-Mukai transform} with kernel $E$ is given by $\Phi_{E}(\bullet)=p_{1\ast}(E\otimes p^{\ast}_{2}(\bullet))$. Let us denote by $I_{i}: \mathcal{A}_{i}\hookrightarrow D^{b}(X)$ the embedding functor from the semi-orthogonal component $\mathcal{A}_i$ to $D^{b}(X)$ and use $I^{\ast}_{i}:D^{b}(X)\to \mathcal{A}_i$ to denote the projection functor. A functor $\Phi: \mathcal{A}_{i}\rightarrow \mathcal{A}_{i}$ is  said to be \emph{of Fourier-Mukai type} if there exists $E\in D^{b}(X\times X)$ such that the following diagram is commutative.
      $$\xymatrix@C=2.5cm{D^{b}(X)\ar[d]^{I^{\ast}_{i}}\ar[r]^{\Phi_{E}}&D^{b}(X)\\
      \mathcal{A}_{i}\ar[r]^{\Phi}&\mathcal{A}_{i}\ar[u]^{I_{i}}}$$
\begin{proposition}\label{FourierMukai}
The following statements are equivalent.
\begin{itemize}
\item The Fourier-Mukai functor $\Phi_{E}$ factors through a semi-orthogonal component $\mathcal{A}_{i}$.
	$$\xymatrix@C=2.5cm{D^{b}(X)\ar[d]^{I^{\ast}_{i}}\ar[r]^{\Phi_{E}}&D^{b}(X)\\
	\mathcal{A}_{i}\ar[r]^{\Phi'_{E}}&\mathcal{A}_{i}\ar[u]^{I_{i}}}$$
\item The kernel $E\in \mathcal{A}_{i}\boxtimes \mathcal{B}^{\vee}_{i}$. 
\end{itemize}
\end{proposition}
\begin{proof}
Firstly, the condition that the functor $\Phi_{E}$ maps $D^{b}(X)$ to $\mathcal{A}_{i}$ is equivalent to the following:
\begin{align*}
&\Hom(p^{\ast}_{1}\mathcal{A}_{j}\boxtimes D^{b}(X), E)\cong \Hom(\mathcal{A}_{j},p_{1\ast}(E\otimes p^{\ast}_{2}D^{b}(X)))=0,\quad j> i \\
&\Hom(E, p_{1}^{!}\mathcal{A}_{s}\boxtimes D^{b}(X))\cong \Hom(p_{1\ast}(E\otimes p^{\ast}_{2}D^{b}(X)), \mathcal{A}_{s})=0,\quad s<i. 
\end{align*}
Since $p^{!}_{1}(\bullet)\cong p^{\ast}_{1}(\bullet)\otimes \omega_{X\times X}\otimes p^{\ast}_{1}\omega^{-1}_{X}[\dim X]$, we have $p^{!}_{1}\mathcal{A}_{s}\boxtimes D^{b}(X)= \mathcal{A}_{s}\boxtimes D^{b}(X)$. Therefore $\Phi_{E}$ mapping $D^{b}(X)$ to $\mathcal{A}_{i}$ is equivalent to $E\in \mathcal{A}_{i}\boxtimes D^{b}(X)$. Secondly, the statement that $\Phi_{E}$ maps $\mathcal{A}_{t}$ to $0$ for $t\neq i$ is equivalent to the following:
\begin{align*}
\Hom(E^{\vee}, D^{b}(X)\boxtimes \mathcal{A}_{t})\cong \Hom(D^{b}(X),p_{1\ast}(E\otimes p^{\ast}_{2}\mathcal{A}_{t}))=0,\quad t\neq i
\end{align*}
which is further equivalent to $E\in D^{b}(X)\boxtimes \mathcal{B}^{\vee}_{i}$.  
\end{proof}

\begin{remark}\label{rmk_dg_cat}
	Let $\operatorname{dg-cat}$ denote the category of $dg$ categories. Write $\operatorname{Hqe}(\operatorname{dg-cat})$ as the localized dg-cat with respect to quasi-equivalences of dg categories. Let $\operatorname{L}_{parf}(X)$ be a dg enhancement of $D^{b}(X)$ whose objects are injective complexes, and let $\mathcal{A}_{idg}\subset\operatorname{L}_{parf}(X)$ be a sub-dg category whose objects are in $\mathcal{A}_{i}\subset D^{b}(X)$. According to \cite[Theorem 8.15]{Toenhomotopic2007}, it holds in $\operatorname{Hqe}(\operatorname{dg-cat})$ that
	$$RHom(\operatorname{L}_{parf}(X), \operatorname{L}_{parf}(X))\cong \operatorname{L}_{parf}(X\times X)$$
where $RHom(-,-)$ is the internal. In particular, the projection $j^{\ast}_{i}E$ with $E\in D^{b}(X\times X)$ can be regarded as a quasi-functor from $\mathcal{A}_{idg}$ to $\mathcal{A}_{idg}$.
\end{remark}
 
\begin{proposition}\label{categoricalFourierMukai}
	Let $\Phi_{1}$ and $\Phi_{2}$ be Fourier-Mukai functors of $\mathcal{A}_{i}$ with kernels $E_{1}$ and $E_{2}$ respectively. Write $D=[RHom(\mathcal{A}_{i},\mathcal{A}_{i})]$. Then 
 $$\Hom_{D}(\Phi_{1},\Phi_{2})\cong \Hom_{D^{b}(X\times X)}(E_{1},E_{2}).$$ 
\end{proposition}
\begin{proof}
    It holds that 
	$$\Hom_{[RHom(\operatorname{L}_{parf}(X),\operatorname{L}_{parf}(X))]}(\Phi_{E_{1}},\Phi_{E_{2}})\cong \Hom_{D}(\Phi_{1},\Phi_{2}).$$ 
By \cite[Theorem 8.15]{Toenhomotopic2007}, we also have
	$$\Hom_{[RHom(\operatorname{L}_{parf}(X), \operatorname{L}_{parf}(X))]}(\Phi_{E_{1}},\Phi_{E_{2}})\cong \Hom_{D^{b}(X\times X)}(E_{1},E_{2}).$$
The proof is then completed.
\end{proof}

\subsubsection{Hochschild (co)homology as Fourier-Mukai kernels}
	Notation as in Subsection \ref{section_FM_functors}. In particular, let $P_{i}=j'^{\ast}_{i}\Delta_{\ast}\mathcal{O}_{X}\cong j^{\ast}\Delta_{\ast}\mathcal{O}_{X}$. In this subsection, we prove the following proposition which allows us to compute Hochschild homology using Fourier-Mukai kernels. (The result holds for any semi-orthogonal component $\mathcal{A}_i$; we state the proposition only for $\mathcal{A}_1$ because later in Section \ref{subsection_equivariant_inf_categorical_Torelli} all the other semi-orthogonal components are exceptional and $\mathcal{A}_1$ denotes the only geometrically meaningful residual subcategory.)
\begin{proposition}\label{Hochschildhomology}
Let $F\in D^{b}(X\times X)$. Then the following holds. 
	$$\Hom(j^{\ast}_{1}F, j^{\ast}_{1}\Delta_{\ast}\mathcal{O}_{X})\cong  \Hom(j'^{\ast}_{1}F, j'^{\ast}_{1}\Delta_{\ast}\mathcal{O}_{X}).$$
In particular,  
$$\HH_{\bullet}(\mathcal{A}_{1})\cong\Hom(P_{1}\circ S^{-1}_{X}, P_{1}[\bullet])\cong \Hom(j_{1}^{\ast}S^{-1}_{X}, j_{1}^{\ast}\Delta_{\ast}\mathcal{O}_{X}[\bullet]).$$
\end{proposition}

	We divide the proof of Proposition \ref{Hochschildhomology} into the following several lemmas. 
\begin{lemma}\label{convolutionbox}
	Let $E\in \mathcal{A}_{s_{1}}\boxtimes \mathcal{B}^{\vee}_{t_{1}}$ and let $F\in \mathcal{A}_{s_{2}}\boxtimes \mathcal{B}^{\vee}_{t_{2}}$. Then $E\circ F\in \mathcal{A}_{s_{1}}\boxtimes \mathcal{B}^{\vee}_{t_{2}}$. Moreover, $E\circ F=0$ whenever $t_{1}\neq s_{2}$.
\end{lemma}
\begin{proof}
WLOG, We assume that $E=E_{s_{1}}\boxtimes E_{t_{1}}$ and that $F=F_{s_{2}}\boxtimes F_{t_{2}}$ where $E_{s_{1}}\in \mathcal{A}_{s_{1}}$, $E_{t_{1}}\in \mathcal{B}^{\vee}_{t_{1}}$, $F_{s_{2}}\in \mathcal{A}_{s_{2}}$ and $F_{t_{2}}\in \mathcal{B}^{\vee}_{t_{2}}$. Denote by $p'_{i}$ the $i$-th projection from $X\times X\times X$ to $X$ with $i=1,2,3$. 
 
$$\xymatrix{&X\times X\times X\ar[dl]_{p_{12}}\ar[dr]^{p_{23}}&\\
X\times X&&X\times X}$$
Then 
\begin{align*}
    E\circ F=\, &p_{13\ast}(p^{\ast}_{12}(p^{\ast}_{1}E_{s_{1}}\otimes p^{\ast}_{2}E_{t_{1}})\otimes p^{\ast}_{23}(p^{\ast}_{1}F_{s_{2}}\otimes p^{\ast}_{2}F_{t_{2}}))\\
            \cong\,& p_{13\ast}(p'^{\ast}_{1}E_{s_{1}}\otimes p'^{\ast}_{2}E_{t_{1}}\otimes p'_{2}F_{s_{2}}\otimes p'^{\ast}_{3}F_{t_{2}})\\
            \cong\,& p_{13\ast}(p^{\ast}_{13}(p^{\ast}_{1}E_{s_{1}}\otimes p^{\ast}_{2}F_{t_{2}})\otimes p'^{\ast}_{2}(E_{t_{1}}\otimes F_{s_{2}}))\\
            \cong\,& p^{\ast}_{1}E_{s_{1}}\otimes p^{\ast}_{2}F_{t_{2}}\otimes p_{13\ast}p'^{\ast}_{2}(E_{t_{1}}\otimes F_{s_{2}})\\
            \cong\,& p^{\ast}_{1}E_{s_{1}}\otimes p^{\ast}_{2}F_{t_{2}}\otimes R\Gamma(E_{t_{1}}\otimes F_{s_{2}}).
\end{align*}
(For the last isomorphism, consider the following cartesian diagram.
	$$\xymatrix@C=2cm{X\times X\times X\ar[d]_{p_{13}}\ar[r]^(0.6){p'_{2}}& X\ar[d]\\
X\times X\ar[r]&k}$$
Hence $p_{13\ast}p'_{2}(E_{t_{1}\otimes F_{s_{2}}})\cong R\Gamma(E_{t_{1}}\otimes F_{s_{2}})\otimes \mathcal{O}_{X\times X}$.)
Since $R\Gamma(E_{t_{1}}\otimes F_{s_{2}})\cong \mathrm{RHom}(E^{\vee}_{t_{1}}, F_{s_{2}})=0$ for any $t_{1}\neq s_{2}$, one gets $E\circ F=0$ whenever $t_{1}\neq s_{2}$. 
\end{proof}

\begin{lemma}\label{lemma1}
Let $F\in D^{b}(X\times X)$. Recall that $P_{i}=j'^{\ast}_{i}\Delta_{\ast}\mathcal{O}_{X}\cong j^{\ast}\Delta_{\ast}\mathcal{O}_{X}$. Then the following hold.
\begin{enumerate}
\item $j'^{\ast}_{i}F\cong P_{i}\circ F$. 
\item $j^{\ast}_{i}F\cong P_{i}\circ F\circ P_{i}$.  
\end{enumerate} 
\end{lemma}
\begin{proof}
	For the first isomorphism, we have $j'^{\ast}_{i}F\cong \Delta_{\ast}\mathcal{O}_{X}\circ j'^{\ast}_{i}F\cong P_{i}\circ j'^{\ast}_{i}F\cong P_{i}\circ F$. For the second one, by Lemma \ref{convolutionbox} one has $j^{\ast}_{i}F\cong j^{\ast}_{i}F\circ \Delta_{\ast}\mathcal{O}_{X}\cong j^{\ast}_{i}F\circ P_{i}$. It then follows that $j^{\ast}_{i}F\cong \Delta_{\ast}\mathcal{O}_{X}\circ j^{\ast}_{i}F\circ P_{i}\cong P_{i}\circ j^{\ast}_{i}F\circ P_{i}\cong P_{i}\circ F\circ P_{i}$ (the last equality holds because $P_{i}\circ F_{s,t}\circ P_{i}=0$ for the semi-orthogonal factor $F_{s,t}$ of $F$ in $\mathcal{A}_{s}\boxtimes \mathcal{B}^{\vee}_{t}$ with $s\neq i$ or $t\neq i$. 
\end{proof}

\begin{lemma}\label{jrestriction}
	For any $E\in D^{b}(X\times X)$, $\Phi_{j^{\ast}_{1}E}\vert_{\mathcal{A}_{1}}=\Phi_{j'^{\ast}_{1}E}\vert_{\mathcal{A}_{1}}$. Moreover, $j^{\ast}_{1}S^{-1}_{X}$ is the inverse of Serre functor of $\mathcal{A}_{1}$.  
\end{lemma}

\begin{proof}
For any object $E\in D^{b}(X\times X)$ and $A_{1}\in\mathcal{A}_{1}$, there exists a triangle
	$$B\rightarrow j'^{\ast}_{1}E\rightarrow j^{\ast}_{1}E$$
where $B$ is generated by $\{\mathcal{A}_{1}\boxtimes \mathcal{B}^{\vee}_{t}\}_{t\geq 2}$. Since $\Phi_{B}(A_{1})\cong p_{1\ast}(B\otimes p^{\ast}_{2}(A_{1}))\cong \mathcal{A}_{1}\otimes RHom(\mathcal{B}_{t}, A_{1})=0$, we have  
	$$\Phi_{j^{\ast}_{1}(E)}(A_{1})\cong \Phi_{j'^{\ast}_{1}E}(A_{1}).$$
In other words, $\Phi_{j^{\ast}_{1}E}\vert_{\mathcal{A}_{1}}=\Phi_{j'^{\ast}_{1}E}\vert_{\mathcal{A}_{1}}$. By adjunction, $\Phi_{j'^{\ast}_{1}S^{-1}_{X}}(E)$ co-represents the functor $F\mapsto \Hom(E,F)^{\vee}$. Since the Serre functor also represents the functor $F\mapsto \Hom(E,F)^{\vee}$, we conclude that $\Phi_{j^{\ast}_{1}S^{-1}_{X}}$ is the inverse of Serre functor of $\mathcal{A}_{1}$.   
\end{proof}

\begin{proof}[Proof of Proposition \ref{Hochschildhomology}.] 
Let us first prove $\Hom(j^{\ast}_{1}F,j^{\ast}_{1}\Delta_{\ast}\mathcal{O}_{X})\cong \Hom(j'^{\ast}_{1}F,j'^{\ast}_{1}\Delta_{\ast}\mathcal{O}_{X})$. According to Lemma \ref{lemma1}, it suffices to show that $\Hom(P_{1}\circ F\circ P_{j}, P_{1})=0$ for $j\neq 1$. Let $F_{s,t}$ be the semi-orthogonal factor of $F$ in $\mathcal{A}_{s}\boxtimes \mathcal{B}^{\vee}_{t}$. By Lemma \ref{convolutionbox}, $P_{1}\circ F_{s,t}\circ P_{j}\in \mathcal{A}_{1}\boxtimes\mathcal{B}^{\vee}_{j}$. By semi-orthogonality, we get $\Hom(P_{1}\circ F_{s,t}\circ P_{j}, P_{1})=0$ which implies that $\Hom(P_{1}\circ F\circ P_{j}, P_{1})=0$. Next, according to Lemma \ref{jrestriction}, $j^{\ast}_{1}S^{-1}_{X}$ is the inverse Serre functor of $\mathcal{A}_{1}$. We then have
\begin{align*}
    \HH_{\bullet}(\mathcal{A}_{1})=\Hom(S^{-1}_{\mathcal{A}_{1}},\operatorname{Id}_{\mathcal{A}_{1}}[\bullet])\cong\,& \Hom(j^{\ast}_{1}S^{-1}_{X},j^{\ast}_{1}\Delta_{\ast}\mathcal{O}_{X}[\bullet])\\
    \cong\,& \Hom(P_1\circ S^{-1}_{X}\circ P_1, P_1[\bullet])\\
    \cong\, &\Hom(P_1\circ S^{-1}_{X},P_1[\bullet]).
\end{align*} 
\end{proof}

\subsubsection{Compatibility of Hochschild (co)homology}
	Notation as in Subsection \ref{section_FM_functors}. In this subsection, we aim at comparing the algebra structures of the Hochschild cohomology $\HH^{\ast}(X)$ and $\HH^{\ast}(\mathcal{A}_{1})$, and comparing the module structures on the Hochschild homology $\HH_{\ast}(X)$ and $\HH_{\ast}(\mathcal{A}_{1})$, using Fourier-Mukai transforms. From now on, we always assume that  $\Phi_{H}: D^{b}(X)\rightarrow D^{b}(X)$ is a Fourier-Mukai transform which preserves each semi-orthogonal component $\mathcal{A}_{i}$, and we write $\Phi_{iH}$ to denote the induced functor on $\mathcal{A}_{i}$. We refer the reader to \cite[\S4.3]{Caldararu2010} and \cite[\S4.1]{perryhochschildgroup20211}\footnote{Although the author only wrote down the definition for pushforward for Hochschild cohomology, it is easy to see pushforward of Hochschild homology can be defined similarly.} for the definitions of the morphisms $\Phi_{H\ast}: \HH_{\ast}(X)\rightarrow \HH_{\ast}(X)$ and $\Phi_{iH\ast}:\HH_{\ast}(\mathcal{A}_{i})\rightarrow \HH_{\ast}(\mathcal{A}_{i})$ via Fourier-Mukai kernels. 

\begin{lemma}\label{commuteprojection} 
	Notation as above. Suppose that the Fourier-Mukai transform $\Phi_{H}: D^{b}(X)\rightarrow D^{b}(X)$ preserves each semi-orthogonal component $\mathcal{A}_{i}$. Then for any $i$ the kernel $H$ satisfies the following properties:
\begin{enumerate}
\item $j'^{\ast}_{i}H=j^{\ast}_{i}H$;
\item $E\circ H\cong E\circ j'^{\ast}_{i}H\cong E\circ j^{\ast}_{i}H$ where $E\in \mathcal{A}_{i}\boxtimes \mathcal{B}^{\vee}_{i}$;
\item $H\circ  P_{i}\cong P_{i}\circ H$.
\end{enumerate}
Moreover, suppose that both of the Fourier-Mukai transforms $\Phi_{H_{1}}$ and $\Phi_{H_{2}}$ preserve all the semi-orthogonal components $\mathcal{A}_{i}$. Also let $F\in D^{b}(X\times X)$. Then it holds that $$j^{\ast}_{i}(H_{1}\circ F\circ H_{2})\cong j^{\ast}_{i}H_{1}\circ j^{\ast}_{i}F\circ j^{\ast}_{i}H_{2}$$ which is functorial with respect to the factor $F$. 
\end{lemma}
\begin{proof}
For (1), according to Lemma \ref{intersection} and Lemma \ref{lemma1}, it suffices to prove that $P_{i}\circ H\in _{X}\mathcal{B}^{\vee}_{i}$. Since the composition of Fourier-Mukai functors $\Phi_{P_{i}}\circ \Phi_{H}$ maps $\mathcal{A}_{i}$ to $\mathcal{A}_{i}$, and maps $\mathcal{A}_{j}$ to zero for $j\neq i$, the same argument as that in the proof of Proposition \ref{FourierMukai} shows that $P_{i}\circ H\in _{X}\mathcal{B}^{\vee}_{i}$. Part (2) follows from Lemma \ref{convolutionbox}. For (3), by (1), (2) and Lemma \ref{convolutionbox},
	$$P_{i}\circ H\cong P_{i}\circ j^{\ast}_{i}H\cong \Delta_{\ast}\mathcal{O}_{X}\circ j^{\ast}_{i}H\cong j^{\ast}_{i}H\circ \Delta_{\ast}\mathcal{O}_{X}\cong j^{\ast}_{i}H\circ P_{i}\cong H\circ P_{i}.$$  
For the last statement, by (3) and Lemma \ref{lemma1},
 \begin{align*}
     j^{\ast}_{i}(H_{1}\circ F\circ H_{2})\cong  
     P_{i}\circ H_{1}\circ F\circ H_{2}\circ P_{i}
     &\cong
     P_{i}\circ H_{1}\circ P_{i}\circ F\circ P_{i}\circ H_{2}\circ P_{i}\\
      &=j^{\ast}_{i}H_{1}\circ j^{\ast}_{i}F\circ j^{\ast}_{i}H_{2}
\end{align*}
which is functorial with respect to the factor $F$. 
\end{proof}

\begin{remark}
	Notation as in Remark \ref{rmk_dg_cat}. There exists an isomorphism between the invariants of hochschild homology
$$\HH_{\ast}(X)^{\Phi_{H}\ast}\cong \bigoplus_{i}\HH_{\ast}(\mathcal{A}_{i})^{\Phi_{iH\ast}}.$$
Indeed, the embedding $\mathcal{A}_{idg}\subset \operatorname{L}_{parf}(X)$ induces the following isomorphism of Hochschild homology 
	$$\HH_{\ast}(\operatorname{L}_{parf}(X))\cong \bigoplus\HH_{\ast}(\mathcal{A}_{idg})$$
which is compatible with the quasi-functors $\Phi_{H\ast}: \operatorname{L}_{parf}(X)\rightarrow \operatorname{L}_{parf}(X)$ and $\Phi_{iH\ast}: \mathcal{A}_{idg}\rightarrow \mathcal{A}_{idg}$. Namely we have the following commutative diagram
$$\xymatrix{\HH_{\ast}(\operatorname{L}_{parf}(X))\ar[r]^{\simeq}\ar[d]_{\Phi'_{H\ast}}&\bigoplus\HH_{\ast}(\mathcal{A}_{idg})\ar[d]^{(\Phi'_{1H\ast},\Phi'_{2H\ast},\cdots,\Phi'_{mH\ast})}\\
\HH_{\ast}(\operatorname{L}_{parf}(X))\ar[r]^{\simeq}& \bigoplus\HH_{\ast}(\mathcal{A}_{idg})}$$
where $\Phi'_{H\ast}$ and $\Phi'_{iH\ast}$ are morphisms of Hochschild homology as defined in \cite[Definition 2.1]{Polichuklefchetzdg}. According to \cite[Corollary 2.11]{Polichuklefchetzdg}, we have $\Phi'_{H\ast}=\Phi_{H\ast}$ and $\Phi'_{iH\ast}=\Phi_{iH\ast}$ as morphisms of Hochschild homology. The above result then follows.
\end{remark}

	In what follows, we assume in addition that the Fourier-Mukai transform $\Phi_{H}:D^b(X)\to D^b(X)$ is an auto-equivalence preserving all the semi-orthogonal components. Denote the inverse auto-equivalence by $\Phi_{H^{-1}}$. Now let us make the following observations. 
\begin{itemize}
\item The auto-equivalence $\Phi$ induces an automorphism of Hochschild homology $\Phi_{H\ast}: \mathrm{HH}_{\ast}(X)\rightarrow \mathrm{HH}_{\ast}(X)$. More precisely, for any $a\in \mathrm{HH}_{m}(X)$, $\Phi_{H\ast}(a)$ is given by the following composition. 
\begin{equation}\label{eqn_Phi_*H}
\xymatrix@C+1pc{S^{-1}_{X}[-m]\ar[r]&H\circ S^{-1}_{X}[-m]\circ H^{-1}\ar[r]^{Id\circ a\circ Id}&H\circ \Delta_{\ast}\mathcal{O}_{X}\circ H^{-1}\ar[r]&\Delta_{\ast}\mathcal{O}_{X}}
\end{equation}
\item The projection $j^{\ast}_{1}: D^{b}(X\times X)\rightarrow \mathcal{A}_{1}\boxtimes \mathcal{B}^{\vee}_{1}$ gives rise to a morphism $j^{\ast}_{1}:\HH_{\ast}(X)\rightarrow \HH_{\ast}(\mathcal{A}_{1})$ (as a consequence of  Corollary \ref{Hochschildhomology}).
\item Let $\Phi_{1H}: \mathcal{A}_{1}\rightarrow \mathcal{A}_{1}$ be the auto-equivalence of $\mathcal{A}_1$ induced by $\Phi_{H}$. Then there is an automorphism $\Phi_{1H\ast}:\mathrm{HH}_{\ast}(\mathcal{A}_{1})\rightarrow \mathrm{HH}_{\ast}(\mathcal{A}_{1})$ defined as the action on natural transformations of functors (see \cite[\S4.3]{Caldararu2010} and \cite[\S2.1]{Polichuklefchetzdg}).
\end{itemize}

\begin{proposition}\label{prop_comm_hochschild_homology}
	There exists a commutative diagram as follows.
	$$\xymatrix{\mathrm{HH}_{\ast}(X)\ar[r]^{\Phi_{H\ast}}\ar[d]_{j^{\ast}_{1}}&\mathrm{HH}_{\ast}(X)\ar[d]^{j^{\ast}_{1}}\\
   \mathrm{HH}_{\ast}(\mathcal{A}_{1})\ar[r]^{\Phi_{1H\ast}}&\mathrm{HH}_{\ast}(\mathcal{A}_{1})}$$
\end{proposition}
\begin{proof}
	Applying the functor $j^{\ast}_{1}: D^{b}(X\times X)\rightarrow \mathcal{A}_{1}\boxtimes \mathcal{B}^{\vee}_{1}$ to the sequence \eqref{eqn_Phi_*H}, we obtain 
	$$\xymatrix@C+1pc{j^{\ast}_{1}S^{-1}_{X}[-m]\ar[r]&j^{\ast}_{1}(H\circ S^{-1}_{X}[-m]\circ H^{-1}\ar[r]^{Id\circ a\circ Id})\ar[r]&j^{\ast}_{1}(H\circ \Delta_{\ast}\mathcal{O}_{X}\circ H^{-1}\ar[r])\ar[r]&j^{\ast}_{1}(\Delta_{\ast}\mathcal{O}_{X})}.$$
According to Lemma \ref{commuteprojection} (4), we have natural isomorphisms $j^{\ast}_{1}(H\circ \Delta_{\ast}\mathcal{O}_{X}\circ H^{-1})\cong j^{\ast}_{1}H\circ P_{1}\circ j^{\ast}_{1}H^{-1}$ and $j^{\ast}_{1}(H\circ S^{-1}_{X}[-m]\circ H^{-1})\cong j^{\ast}_{1}H\circ j^{\ast}_{1}S^{-1}_{X}[-m]\circ j^{\ast}_{1}H^{-1}$. By Proposition \ref{categoricalFourierMukai} and Lemma \ref{jrestriction}, $j^{\ast}_{1}S^{-1}_{X}$ is the inverse Serre functor of $\mathcal{A}_{1}$, $P_{1}$ is the identity functor, and $j^{\ast}_{1}H$ is the functor corresponding to $\Phi_{1H}$. The proposition then follows.
\end{proof}	

\begin{remark}
	Note that $P_{1}\circ H\cong H\circ P_{1}$. By \cite[Lemma 6.7]{kuznetsov2009hochschild}, we have a morphism $\phi_{H}: \mathrm{HH}_{\ast}(\mathcal{A}_{1})\rightarrow \mathrm{HH}_{\ast}(\mathcal{A}_{1})$. Moreover, the morphism $\phi_H$ constructed by Kuznetsov coincides with the above $\Phi_{1H\ast}$.	
\end{remark}

 	Similarly for Hochschild cohomology, one also has the induced morphisms $\Phi_{H\ast}$ and $j^{\ast}_{1}$.
\begin{proposition}\label{Hochschildfunctorial}
 	We have the following commutative diagram.
	$$\xymatrix{\mathrm{HH}^{\ast}(X)\ar[r]^{\Phi_{H\ast}}\ar[d]_{j^{\ast}_{1}}&\mathrm{HH}^{\ast}(X)\ar[d]^{j^{\ast}_{1}}\\
   \mathrm{HH}^{\ast}(\mathcal{A}_{1})\ar[r]^{\Phi_{1H\ast}}&\mathrm{HH}^{\ast}(\mathcal{A}_{1})}$$
\end{proposition}

    Putting the above results together, we obtain the following compactibility theoerm.
 \begin{theorem}\label{compatibilitymain}
	The morphisms $\Phi_{H\ast}$, $\Phi_{1H\ast}$ and $j^{\ast}_{1}$ in Propositions \ref{prop_comm_hochschild_homology} and \ref{Hochschildfunctorial} are compatible with respect to the algebra structure of Hochschild cohomology and the module structure of Hochschild homology over Hochschild cohomology respectively. Namely, the following diagrams are commutative.
  \[\begin{tikzcd}
	&& \mathrm{HH}^{\ast}(\cA_1)\times\mathrm{HH}^{\ast}(\cA_1) & \mathrm{HH}^{\ast}(\cA_1)\times\mathrm{HH}^{\ast}(\cA_1) \\
	\mathrm{HH}^{\ast}(X)\times\mathrm{HH}^{\ast}(X) & \mathrm{HH}^{\ast}(X)\times\mathrm{HH}^{\ast}(X) & \mathrm{HH}^{\ast}(\cA_1) & \mathrm{HH}^{\ast}(\cA_1) \\
	\mathrm{HH}^{\ast}(X) & \mathrm{HH}^{\ast}(X)
	\arrow[from=2-1, to=1-3, "j_1^*"]
	\arrow[from=2-2, to=1-4, "j_1^*"]
	\arrow[dashed, from=3-1, to=2-3, "j_1^*"]
	\arrow[from=3-2, to=2-4, "j_1^*"]
	\arrow[from=2-1, to=3-1, "\cup"]
	\arrow[from=2-2, to=3-2, "\cup"]
	\arrow[from=2-1, to=2-2, "\Phi_{H_*}"]
	\arrow[from=3-1, to=3-2, "\Phi_{H_*}"]
	\arrow[from=1-3, to=1-4, "\Phi_{1H_*}"]
	\arrow[from=2-3, to=2-4, "\Phi_{1H_*}"]
	\arrow[from=1-4, to=2-4, "\cup"]
	\arrow[dashed, from=1-3, to=2-3, "\cup"]
\end{tikzcd}\]

\[\begin{tikzcd}
	&& \mathrm{HH}^{\ast}(\cA_1)\times\mathrm{HH}_*(\cA_1) & \mathrm{HH}^{\ast}(\cA_1)\times\mathrm{HH}_*(\cA_1) \\
	\mathrm{HH}^{\ast}(X)\times\mathrm{HH}_*(X) & \mathrm{HH}^{\ast}(X)\times\mathrm{HH}_*(X) & \mathrm{HH}_*(\cA_1) & \mathrm{HH}_*(\cA_1) \\
	\mathrm{HH}_*(X) & \mathrm{HH}_*(X)
	\arrow[from=2-1, to=1-3, "j_1^*"]
	\arrow[from=2-2, to=1-4, "j_1^*"]
	\arrow[dashed, from=3-1, to=2-3, "j_1^*"]
	\arrow[from=3-2, to=2-4, "j_1^*"]
	\arrow[from=2-1, to=3-1, "\cap"]
	\arrow[from=2-2, to=3-2, "\cap"]
	\arrow[from=2-1, to=2-2, "\Phi_{H_*}"]
	\arrow[from=3-1, to=3-2, "\Phi_{H_*}"]
	\arrow[from=1-3, to=1-4, "\Phi_{1H_*}"]
	\arrow[from=2-3, to=2-4, "\Phi_{1H_*}"]
	\arrow[from=1-4, to=2-4, "\cap"]
	\arrow[dashed, from=1-3, to=2-3, "\cap"]
\end{tikzcd}\]
\end{theorem}
\begin{proof}
	We will only check the commutativity of the first diagram. The second one follows in a similar manner. After applying Proposition \ref{Hochschildfunctorial}, it remains to verify the following: let $a\in HH^{t}(X)$ and $b\in HH^{s}(X)$, then the product $ab\in HH^{s+t}(X)$ is given by the composition  
$$\xymatrix{\Delta_{\ast}\mathcal{O}_{X}\ar[r]^{b}& \Delta_{\ast}\mathcal{O}_{X}\ar[r]^{a}[s]&\Delta_{\ast}\mathcal{O}_{X}[s+t]}.$$
Since $j^{\ast}_{1}\circ \cup(a,b)$ is the composition 
$$\xymatrix{j^{\ast}_{1}\Delta_{\ast}\mathcal{O}_{X}\ar[r]^{j^{\ast}_{1}b}& j^{\ast}_{1}\Delta_{\ast}\mathcal{O}_{X}\ar[r]^{j^{\ast}_{1}a}[s]&j^{\ast}_{1}\Delta_{\ast}\mathcal{O}_{X}[s+t]},$$
one gets $j^{\ast}_{1}\circ \cup (a,b)=\cup(j^{\ast}_{1}a,j^{\ast}_{1}b)$ which completes the proof.
\end{proof}

\begin{corollary}\label{commutativeequivariant}
	The following diagram is commutative.
	$$\xymatrix{\mathrm{HH}^{m}(X)^{\Phi_{H\ast}}\times \mathrm{HH}_{n}(X)^{\Phi_{H\ast}}\ar[r]^(0.6){\cap}\ar[d]_{j^{\ast}_{1}}&\mathrm{HH}_{n+m}(X)^{\Phi_{H\ast}}\ar[d]^{j^{\ast}_{1}}\\
\mathrm{HH}^{m}(\mathcal{A}_{1})^{\Phi_{1H\ast}}\times \mathrm{HH}_{n}(\mathcal{A}_{1})^{\Phi_{1H\ast}}\ar[r]^(0.6){\cap}&\mathrm{HH}_{n+m}(\mathcal{A}_{1})^{\Phi_{1H\ast}}}$$
\end{corollary}

\subsection{Equivariant HKR}\label{sec_equivariant_HKR}
	Let $X$ be a smooth projective variety which admits a biregular automorphism $\tau$. The pullback functor induces an action $\tau^*$ on Hochschild (co)homology. There is also a natural action $\tau^*$ on the harmonic structure of the singular cohomology
	$$(\mathrm{HT}^{\bullet}(X)=\bigoplus_{p+q=\bullet}H^{p}(X,\wedge^{q}T_{X}), \ \mathrm{H\Omega}_{\bullet}(X)=\bigoplus_{p-q=\bullet}H^{p}(X,\Omega^{q}_{X})).$$ 
By \cite[Theorem 1.2]{macri2009infinitesimal}, the  $\operatorname{HKR}$ isomorphism for Hochschild homology is compatible with the actions $\tau^*$; more generally, similar compatibility result holds for actions induced by Fourier-Mukai functors. In this subsection, we show that the twisted $\operatorname{HKR}$ isomorphism for Hochschild cohomology is also compatible with the actions induced by $\tau$. 

	We start by recalling the definitions of Hochschild homology and cohomology. The Hochschild cohomology and homology of $X$ are defined respectively as
	$$\mathrm{HH}^{\ast}(X)=\mathrm{Ext}^{\ast}(\Delta_{\ast}\mathcal{O}_{X},\Delta_{\ast}\mathcal{O}_{X});\,\,\,\,\,\, 
\mathrm{HH}_{\ast}(X)=\mathrm{Ext}^{\ast}(\mathcal{O}_{X\times X}, \Delta_{\ast}\mathcal{O}_{X}\otimes^{L}\Delta_{\ast}\mathcal{O}_{X}).$$
Note that $\mathrm{HH}_{\ast}(X)$ is a graded module over the graded algebra $\mathrm{HH}^{\ast}(X)$.
\begin{proposition}\label{IK}{(\cite[Theorem 1.4]{calaque2012cualduararu})}
We have the following twisted $\operatorname{HKR}$ isomorphisms 
\begin{equation*}
\operatorname{IK}: \mathrm{HH}^{\ast}(X)\cong \bigoplus_{p+q=\ast}H^{p}(X,\wedge^{q}T_{X})
\end{equation*}
and
\begin{equation*}
 \operatorname{IK}: \mathrm{HH}_{\ast}(X)\cong \bigoplus_{p-q=\ast}H^{p}(X,\Omega^{q}_{X}).   
\end{equation*}
The isomorphism $\operatorname{IK}$ for Hochschild cohomology is an isomorphism of algebras, where the right hand side is an algebra of poly-vector fields. The isomorphism $\operatorname{IK}$ for Hochschild homology is an isomorphism of graded modules, where the differential forms naturally form a module over the algebra of poly-vector fields.
\end{proposition}
    
	Define $\tau': X\rightarrow X\times X$ as $\tau'(x)=(x,\tau x)$ and let
$(\tau,Id):X\times X\rightarrow X\times X$ be given by $\tau (x,y)=(\tau x,y)$
\begin{lemma}
The pullback functor $\tau^{\ast}:D^{b}(X)\rightarrow D^{b}(X)$ is isomorphic to the Fourier-Mukai functors $\Phi_{\tau'_{\ast}\mathcal{O}_{X}}\cong \Phi_{(\tau,Id)^{\ast}(\Delta_{\ast}\mathcal{O}_{X})}$.
\end{lemma}
\begin{proof}
On one hand, we have the following isomorphisms:
\begin{equation*}
  \Phi_{\tau'_{\ast}\mathcal{O}_{X}}(\bullet)\cong p_{1\ast}(\tau'_{\ast}(\mathcal{O}_{X})\otimes p_{2}^{\ast}(\bullet))
  \cong p_{1\ast}\tau'_{\ast}(\mathcal{O}_{X}\otimes \tau'^{\ast}p^{\ast}_{2}(\bullet)) 
  \cong \tau^{\ast}(\bullet).
\end{equation*}
On the other hand, consider the Cartesian diagram below.
$$\xymatrix{X\ar[r]^{\tau}\ar[d]^{\tau'}&X\ar[d]^{\Delta}\\
X\times X\ar[r]^{(\tau,Id)}&X\times X}$$   
Thus we get $\tau'_{\ast}\mathcal{O}_{X}\cong\tau'_{\ast}\tau^{\ast}\mathcal{O}_{X}\cong (\tau,Id)^{\ast}\Delta_{\ast}\mathcal{O}_{X}$. The lemma then follows.
\end{proof}

\begin{theorem}\label{IKequivariant}
	Let $X$ be a smooth projective variety with an automorphism $\tau$. Then the twisted $\operatorname{HKR}$ isomorphism for Hochschild cohomology is compatible with the actions induced by $\tau$. In other words, the following diagram is commutative.
	$$\xymatrix{\mathrm{HH}^{\ast}(X)\ar[r]^{\tau^{\ast}}\ar[d]^{IK}&\mathrm{HH}^{\ast}(X)\ar[d]^{IK}\\
\bigoplus_{p+q=\ast}H^{p}(X,\wedge^{q}T_{X})\ar[r]^{\tau^{\ast}}& \bigoplus_{p+q=\ast}H^{p}(X,\wedge^{q}T_{X})}$$
\end{theorem}  
\begin{proof}
	Let $a\in \mathrm{HH}^{m}(X)=\Hom(\Delta_{\ast}\mathcal{O}_{X},\Delta_{\ast}\mathcal{O}_{X}[m])$. Write $a'=\operatorname{IK}(a)$, we need to verify that $\tau^{\ast}(a)'=\tau^{\ast}(a')$. By definition, the following diagram is commutative.
	$$\xymatrix@C=3cm{(\tau, Id)^{\ast}(\tau^{-1}, Id)^{\ast}\Delta_{\ast}\mathcal{O}_{X}\ar[r]^{(\tau, Id)^{\ast}(\tau^{-1}, Id)^{\ast}a}&(\tau, Id)^{\ast}(\tau^{-1}, Id)^{\ast}\Delta_{\ast}\mathcal{O}_{X}[m]\ar[d]^{\cong}\\
\Delta_{\ast}\mathcal{O}_{X}\ar[r]^{\tau^{\ast}(a)}\ar[u]^{\cong}&\Delta_{\ast}\mathcal{O}_{X}[m]}.$$
Using adjunction, we also get a commutative diagram as follows. Note that $\Delta^{!}(\bullet)\cong \Delta^{\ast}(\bullet)\otimes \omega^{-1}_{X}[-n]$ where $n$ is the dimension of $X$.
	$$\xymatrix@C=3.5cm{\Delta^{!}(\tau, Id)^{\ast}(\tau^{-1}, Id)^{\ast}\Delta_{\ast}\mathcal{O}_{X}\ar[r]^{(\tau, Id)^{\ast}(\tau^{-1}, Id)^{\ast}a}&\Delta^{!}(\tau, Id)^{\ast}(\tau^{-1}, Id)^{\ast}\Delta_{\ast}\mathcal{O}_{X}[m]\ar[d]\\
\mathcal{O}_{X}\ar[r]^{\tau^{\ast}(a)}\ar[u]&\Delta^{!}\Delta_{\ast}\mathcal{O}_{X}[m]}.$$
Define $\tau'': X\rightarrow X\times X$ as $\tau''x=(\tau x, x)$. From the following we deduce that $\Delta^{\ast}(\tau, Id)^{\ast}=\tau''^{\ast}$.
	$$\xymatrix{&X\ar[d]^{\Delta}\ar[dl]_{\tau''}\\
X\times X&X\times X\ar[l]^{(\tau, Id)}}$$
Since $(\tau^{-1}, Id)^{\ast}\Delta_{\ast}\mathcal{O}_{X}\cong \mathcal{O}_{\Gamma}$ where $\Gamma$ is the subvariety $(x,\tau^{-1}x)\vert_{x\in X}=(\tau x,x)\vert_{x\in X}$, we get $(\tau^{-1}, Id)^{\ast}\Delta_{\ast}\mathcal{O}_{X}\cong (Id, \tau)^{\ast}\Delta_{\ast}\mathcal{O}_{X}$. Also consider the following diagram.
	$$\xymatrix{X\ar[r]^{\tau}\ar[d]^{\tau''}&X\ar[d]^{\Delta}\\
X\times X\ar[r]^{(Id, \tau)}&X\times X}$$
Hence, $\tau''(Id, \tau)^{\ast}=\tau^{\ast}\Delta^{\ast}$. Therefore, we obtain the following commutative diagram. 
	$$\xymatrix@C=3.5cm{\tau^{\ast}\Delta^{\ast}\Delta_{\ast}\mathcal{O}_{X}\otimes \omega^{-1}_{X}[-n]\ar[r]^{\tau^{\ast}\Delta^{\ast}\circ a\circ id}&\tau^{\ast}\Delta^{\ast}\Delta_{\ast}\mathcal{O}_{X}[m]\otimes \omega^{-1}_{X}[-n]\ar[d]\\
\mathcal{O}_{X}\ar[r]^{\tau^{\ast}a}\ar[u]&\Delta^{\ast}\Delta_{\ast}\mathcal{O}_{X}[m]\otimes \omega^{-1}_{X}[-n]}.$$
According to Proposition \ref{IK}, $\tau^{\ast}(a)'=\tau^{\ast}1\tau^{\ast}(a')\tau^{\ast}1=\tau^{\ast}a'$ which completes the proof.
\end{proof}
\begin{remark}
	Together with the results in \cite{macri2009infinitesimal}, a similar argument shows that $IK$ is equivariant with respect to the algebra and module structures. 
\end{remark}

\subsection{Infinitesimal categorical Torelli theorem: An equivariant version}\label{subsection_equivariant_inf_categorical_Torelli}
	In this section, we combine the results in the previous subsections to prove an equivariant version of the infinitesimal categorical Torelli theorem. Let $X$ be a smooth projective variety with an automorphism $\tau$. Assume we have a semi-orthogonal decomposition 
	$$D^{b}(X)=\Big\langle \mathcal{A}_{1},E_{1},E_{2},\cdots,E_{m}\Big\rangle$$
which is preserved by $\tau$. Here $\{E_{i}\}$ is a exceptional collection. According to Corollary \ref{commutativeequivariant} and Theorem \ref{IKequivariant}, we have the following commutative diagram (with notation the same as in Subsection \ref{sec_equivariant_HKR}).
	$$\xymatrix@C=3cm{HT^{2}(X)^{\tau}\times H\Omega_{-1}(X)^{\tau}\ar[r]&H\Omega_{1}(X)^{\tau} \\
   \mathrm{HH}^{2}(X)^{\tau}\times \mathrm{HH}_{-1}(X)^{\tau}\ar[d]^{j^{\ast}_{1}}\ar[r]^(0.6){\cap}\ar[u]^{\operatorname{IK}}&\mathrm{HH}_{1}(X)^{\tau}\ar[d]^{j^{\ast}_{1}}\ar[u]^{\operatorname{IK}}\\
  \mathrm{HH}^{2}(\mathcal{A}_{1})^{\tau}\times \mathrm{HH}_{-1}(\mathcal{A}_{1})^{\tau}\ar[r]^(0.6){\cap}&\mathrm{HH}_{1}(\mathcal{A}_{1})^{\tau}}$$
Since the maps $j^{\ast}_{1}:\mathrm{HH}_1(X)\rightarrow\mathrm{HH}_1(\mathcal{A}_1)$ and $j^{\ast}_1:\mathrm{HH}_{-1}(X)\rightarrow\mathrm{HH}_{-1}(\mathcal{A}_1)$ are isomorphisms and both of them are equivariant with respect to the action $\tau^*$, $j^{\ast}_{1}$ induces isomorphisms $\HH_{-1}(X)^{\tau}\cong \HH_{-1}(\mathcal{A}_{1})^{\tau}$ and $\HH_{1}(X)^{\tau}\cong \HH_{1}(\mathcal{A}_{1})^{\tau}$.
\begin{definition}
	We define $\gamma: \mathrm{HH}^{2}(\mathcal{A}_{1})^{\tau}\rightarrow \Hom(H\Omega_{-1}(X)^{\tau}, H\Omega_{1}(X)^{\tau})$ as follows.
Let $a\in \HH^{2}(\mathcal{A}_{1})^{\tau}$ and $b\in H\Omega_{-1}(X)^{\tau}$, then $\gamma_{a}(b):=IK\circ j^{\ast}_{-1}(a\cap j^{\ast}_{1}\circ \operatorname{IK}^{-1}(b))$. Define $\eta$ as the composition of the inclusion $H^{1}(X,T_{X})\subset HT^{2}(X)$, the isomorphism $\operatorname{IK}$ and the projection $j^{\ast}_{1}$. Define $dp$ as the natural action of vector fields on differential forms.
  \end{definition}

	Combining the commutative diagrams above, we immediately obtain the following theorem.
\begin{theorem}\label{main_theorem_compatibility}
	The following diagram is commutative.
	$$\xymatrix@C=2.5cm{\mathrm{HH}^{2}(\mathcal{A}_{1})^{\tau}\ar[r]^(0.4){\gamma}&\Hom(H\Omega_{-1}(X)^{\tau},H\Omega_{1}(X)^{\tau})\\
H^{1}(X,T_{X})^{\tau}\ar[u]^{\eta}\ar[ur]^(0.4){dp}& }$$
\end{theorem}
%\begin{proof}
%This follows from Corollary \ref{commutativeequivariant} and Theorem \ref{IKequivariant}.
%\end{proof}
 
	To conclude the paper, we consider cubic threefolds with a non-Eckardt involution which have been our main examples. Let $Y$ be a cubic threefold with a non-Ekcardt type involution $\tau$; recall that we have a semi-orthogonal decomposition which is preserved by $\tau$
	$$D^{b}(Y)=\Big\langle \Ku(Y), \mathcal{O}_{Y},\mathcal{O}_{Y}(1)\Big\rangle.$$
Using Theorem \ref{main_theorem_compatibility} we obtain the following corollary. As discussed in Proposition \ref{prop_inf_torelli_nEckardt} (see also \cite[Proposition 3.4]{casalaina2022moduli}), the morphisms $dp$ and hence $\eta$ are injective for a general $(Y,\tau)$. By \cite[Theorem 4.4]{perryhochschildgroup20211}, there is a natural inclusion $\HH^{2}(\Ku(Y))^{\tau}\subset \HH^{2}(\Ku_{\mathbb{Z}_2}(Y))$. Thus Corollary \ref{cor_main_theorem_compatibility} can be regarded as an equivariant version of the infinitesimal categorical Torelli theorem for non-Ekcardt cubic threefolds.    
\begin{corollary}\label{cor_main_theorem_compatibility}
	There exists a commutative diagram for a cubic threefold $Y$ with a non-Eckardt involution $\tau$ as follows.
	$$\xymatrix@C=2.5cm{\mathrm{HH}^{2}(\Ku(Y))^{\tau}\ar[r]^(0.4){\gamma}&\Hom(H\Omega_{-1}(Y)^{\tau},H\Omega_{1}(Y)^{\tau})\\
H^{1}(Y,T_{Y})^{\tau}\ar[u]^{\eta}\ar[ur]^(0.4){dp}& }$$
Moreover, if $(Y,\tau)\in\mathcal{M}_0$ as in Proposition \ref{prop_inf_torelli_nEckardt} (see also \cite[Proposition 3.4]{casalaina2022moduli}), then the map $H^{1}(Y,T_{Y})^{\tau}\stackrel{\eta}{\to} \mathrm{HH}^{2}(\Ku(Y))^{\tau}\subset \HH^{2}(\Ku_{\mathbb{Z}_2}(Y))$ is injective.
\end{corollary}

\appendix\section{An equivariant Categorical Torelli theorem for cubic threefolds}
Let $Y$ be a cubic threefold with a geometric involution $\iota$ (note that $\iota$ could be either of Eckardt type or of non-Eckardt type). In this section we show the equivariant Kuznetsov component $\Ku_{\mathbb{Z}_2}(Y)$, together with the residual action of $\widehat{\mathbb{Z}_2}$, reconstructs the pair $(Y,\iota)$, as suggested by Sasha Kuznetsov and Xiaolei Zhao. 

%it is called Eckardt type if the fix locus of $\iota$ is a point and a cubic surface. This cubic threefold has been studied in detail in \cite{casalaina2021moduli}. After a linear change of coordinates, the equation of $Y$ is given by $f(x_0,\ldots,x_3)+l(x_0,\ldots,x_3)x_4^2=0$, where $f$ is a homogenous polynomial of degree $3$ and $l$ is a linear polynomial. The action of the involution $\iota$ on coordinates is given by 
%$$\iota:[x_0,\ldots,x_3,x_4]\mapsto [x_0,\ldots,x_3,-x_4].$$

%\begin{lemma}\label{lemma_Serre_functor_equivariant_category}
%Let $Y$ be a cubic threefold with an Eckardt type involution $\iota$ and $\Ku_{\mathbb{Z}_2}(Y)$ be the equivariant category of $\Ku(Y)$ under the action of $\iota$. Then the Serre functor $S^{\mathbb{Z}_2}$ of $\Ku_{\mathbb{Z}_2}(Y)$ satisfies $(S^{\mathbb{Z}_2})^3\cong [5]\rho_1$, where $\rho_1$ is the non-trivial irreducible representation of $\mathbb{Z}_2$
%\end{lemma}
%\begin{proof}    
%\end{proof}

%Now, we show the equivariant category $\Ku_{\mathbb{Z}_2}(Y)$ solely reconstruct the isomorphism class of $(Y,\iota)$ whenever $\iota$ is of Eckardt type. 

\begin{theorem}\label{thm_equivariant_cat_Torelli}
Let $\Phi:\Ku_{\mathbb{Z}_2}(Y)\simeq\Ku_{\mathbb{Z}_2}(Y')$ be an equivalence of the equivariant Kuznetsov components. Suppose that $\Phi$ is compatible with the residual actions of $\widehat{\mathbb{Z}_2}$, then $(Y,\iota)\cong (Y',\iota')$. 
\end{theorem}

\begin{proof}
%By Lemma~\ref{lemma_Serre_functor_equivariant_category}, the equivalence $\Phi$ commutes with the residual action by $\widehat{\mathbb{Z}_2}$ on $\Ku_{\mathbb{Z}_2}(Y)$ since $\Phi$ commutes with the Serre functor and shifts and the residual action of $\widehat{\mathbb{Z}_2}$ is given by tensoring with the non-trivial character $\rho_1$. 
Since $\Phi$ commutes with the residual actions of $\widehat{\mathbb{Z}_2}$ on $\Ku_{\mathbb{Z}_2}(Y)$ and $\Ku_{\mathbb{Z}_2}(Y')$, it lifts to an equivalence $$\Psi:\Ku_{\mathbb{Z}_2}(Y)^{\widehat{\mathbb{Z}_2}}\simeq\Ku_{\mathbb{Z}_2}(Y')^{\widehat{\mathbb{Z}_2}}$$ which commutes with the actions of $\widehat{\widehat{\mathbb{Z}_2}}\cong\mathbb{Z}_2$ (cf. the proof of \cite[Lemma 4.9]{bayer2023kuznetsov}). Then by Elagin's reconstruction theorem, it holds that $\Ku_{\mathbb{Z}_2}(Y)^{\widehat{\mathbb{Z}_2}}\simeq\Ku(Y)$, thus we get $$\Psi:\Ku(Y)\simeq\Ku(Y')$$ satisfying $\Psi\circ\iota^*\cong\iota'^*\circ\Psi$, where $\iota^*,\iota'^*$ are the induced auto-equivalences of $\Ku(Y)$ and $\Ku(Y')$ respectively. As a result of \cite[Corollary 6.11]{Li2024Higher}, $\Psi$ is of Fourier-Mukai type. Then by \cite[Theorem 4.23]{liu2024fourier} and \cite[Corollary 8.4]{feyzbakhsh2023new}, there is a unique isomorphism $f:Y\cong Y'$ and integers $m$ and $n$ such that $S_{\Ku(Y)}^n\circ [m]\circ\Psi\cong f_*$. Then we have 
$$S_{\Ku(Y)}^{-n}\circ [-m]\circ f_*\circ\iota^*\cong \iota'^*\circ S_{\Ku(Y)}^{-n}\circ [-m]\circ f_*.$$
It follows that $$f_*\circ\iota^*\cong\iota'^*\circ f_*.$$
But since $\iota^*$ and $\iota'^*$ are both auto-equivalences and involutions, it holds that $\iota^*\circ\iota_*\cong\mathrm{Id}$ and $\iota^*\circ\iota^*\cong\mathrm{Id}$, where $\mathrm{Id}$ is the identical functor on $\Ku(Y)$. Thus $\iota_*\cong\iota^*$. Then we get $f_*\circ\iota_*\cong\iota'_*\circ f_*$ which implies $(f\circ\iota\circ f^{-1})_*\cong\iota'_*$. Using \cite[Remark 8.5]{feyzbakhsh2023new}, we obtain $f\circ\iota\cong\iota'\circ f$ which allows us to conclude that $(Y,\iota)\cong (Y',\iota')$. 
\end{proof}

\begin{remark}
Since $\mathrm{HH}^0(\Ku(Y))\cong \mathbb{C}$ and there is an action of $G=\mathbb{Z}_2$ on $\mathrm{HH}^0(\Ku(Y))\cong \mathbb{C}$, by tensoring with the character of $\mathbb{Z}_2$, denoted by $\kappa$. By \cite[Proposition 5.7, Corollory 5.9]{C19}, the equivariant category $\Ku_{\mathbb{Z}_2}(Y)$ is a fractional Calabi-Yau category of dimension $\frac{10}{6}$, and the Serre functor $S^{\mathbb{Z}_2}$ satisfies $(S^{\mathbb{Z}_2})^3=(\kappa\otimes-)\circ [5]$. If $\kappa$ is the non-trivial character of $\mathbb{Z}_2$, then it is easy to see that the equivalence $\Phi:\Ku_{\mathbb{Z}_2}(Y)\simeq\Ku_{\mathbb{Z}_2}(Y')$ in Theorem \ref{thm_equivariant_cat_Torelli} automatically commutes with residual actions of $\widehat{\mathbb{Z}_2}$. In this case $\Ku_{\mathbb{Z}_2}(Y)$ solely reconstructs the pair $(Y,\iota)$. 
\end{remark}

%\begin{remark}\label{rem_failure_non_Eckardt_type_cubic}
%For non-Eckardt type cubic threefold, the action of $\tau$ on coordinates is given by $\tau: [x_0,\ldots, x_3, x_4,x_5]\mapsto [x_0,\ldots, x_3, -x_4,-x_5]$. In this case $(S^{\mathbb{Z}_2})^3=[5]$. Thus we are not able to lift the equivalence $\Phi:\Ku_{\mathbb{Z}_2}(Y)\simeq\Ku_{\mathbb{Z}_2}(Y')$ to equivalence between $\Ku(Y)$ and $\Ku(Y')$. 
%\end{remark}

%\bibliographystyle{alpha}
%{\small{\bibliography{hochschild}}}
%\nocite{*}

\end{document}